\newif\ifarXiv
\arXivfalse
\newif\ifWP
\WPfalse
\newif\ifFULL
\FULLfalse
\newif\ifCOLOR
\COLORfalse
\newif\ifLATIN
\LATINfalse

\arXivtrue     

\COLORtrue     
\LATINtrue     

\newif\ifnotCOLOR		
\notCOLORtrue
\ifCOLOR\notCOLORfalse\fi

\newif\ifnotLATIN		
\notLATINtrue
\ifLATIN\notLATINfalse\fi

\ifarXiv
  \documentclass{article}
  \usepackage{amsmath,amsthm,amsfonts,latexsym,graphicx,stmaryrd}
  \newcommand{\Extra}[1]{}
\fi

\ifWP
  \documentclass[10pt]{article}
  \input{/Doc/Computing/Latex/kp.txt}
  \usepackage{amsmath,amsthm,amsfonts,latexsym,graphicx,stmaryrd}
  \newcommand{\Extra}[1]{}
\fi

\ifFULL
  \documentclass{article}
  \usepackage{amsmath,amsthm,amsfonts,latexsym,graphicx,stmaryrd,color}
  \newcommand{\Extra}[1]{\blue{#1}}
  
\fi

\ifarXiv
  \title{On-line predictive linear regression}
  \author{Vladimir Vovk, Ilia Nouretdinov, and Alex Gammerman\\
      Computer Learning Research Centre\\
      Department of Computer Science\\
      Royal Holloway, University of London\\
      Egham, Surrey TW20 0EX, UK\\
      \texttt{\{vovk,ilia,alex\}@cs.rhul.ac.uk}}
\fi

\ifWP
  \title{On-line predictive linear regression}
  \author{Vladimir Vovk\and Ilia Nouretdinov\and Alex Gammerman}

  \twodatestrue
  
\fi

\ifFULL
  \title{On-line predictive linear regression}
  \author{Vladimir Vovk, Ilia Nouretdinov, and Alex Gammerman\\
      Computer Learning Research Centre\\
      Department of Computer Science\\
      Royal Holloway, University of London\\
      Egham, Surrey TW20 0EX, UK\\
      \texttt{\{vovk,ilia,alex\}@cs.rhul.ac.uk}}
\fi

\ifFULL\ifCOLOR
  \newcommand{\blue}[1]{\textcolor{blue}{#1}}		
  \newcommand{\makeblue}{\color{blue}}			
\fi\fi
\ifFULL\ifnotCOLOR
  \newcommand{\blue}[1]{\textsl{#1}}
  \newcommand{\makeblue}{\sl}
\fi\fi

\ifLATIN
  \newcommand{\Takeuchi}{takeuchi:1975latin}
\fi
\ifnotLATIN
  \input{OT2enc.def}

  \usepackage{CJK}

  \newcommand{\Takeuchi}{takeuchi:1975}
\fi

\newlength{\picturewidth}
\allowdisplaybreaks[3]

\newcommand{\Vladimir}{Vladimir }
\newcommand{\Alex}{Alex }
\newcommand{\DOT}{.}
\newcommand{\zzrelax}[1]{}

\DeclareMathAlphabet{\mathbfit}{OT1}{cmr}{bx}{it}

\newcommand{\st}{\mathrel{:}}
\newcommand{\given}{\mathrel{|}}

\newcommand{\Zeta}{Z}			
\newcommand{\NM}{\textit{NM}}		

\newcommand{\bbbr}{\mathbb{R}}		
\newcommand{\III}{\mathbb{I}}		
\newcommand{\bbbp}{\mathbb{P}}		
\newcommand{\bbbe}{\mathbb{E}}		
\newcommand{\GGG}{\mathcal{G}}		
\newcommand{\NNN}{\mathcal{N}}		
\newcommand{\PPP}{\mathcal{P}}		

\newcommand{\Prob}{\mathop{\bbbp}\nolimits}
\newcommand{\Expect}{\mathop{\bbbe}\nolimits}

\newcommand{\co}{\mathop{{\rm co}}\nolimits}
\newcommand{\rank}{\mathop{{\rm rank}}\nolimits}
\newcommand{\err}{\mathop{{\rm err}}\nolimits}
\newcommand{\Err}{\mathop{{\rm Err}}\nolimits}

\newcommand{\rnd}{^{\text{{\rm rnd}}}}
\newcommand{\obs}{^{\text{{\rm obs}}}}

\theoremstyle{plain}	
\newtheorem{theorem}{Theorem}
\newtheorem{proposition}{Proposition}
\newtheorem{lemma}{Lemma}
\newtheorem{corollary}{Corollary}

\theoremstyle{remark}
\newtheorem{remark}{Remark}

\theoremstyle{definition}
\newtheorem{definition}{Definition}

\newlength{\IndentI}
\newlength{\IndentII}
\newlength{\IndentIII}
\newlength{\IndentIV}
\newlength{\IndentV}
\newlength{\WidthI}
\newlength{\WidthII}
\newlength{\WidthIII}
\newlength{\WidthIV}
\newlength{\WidthV}
\begin{document}
\maketitle

\begin{abstract}
  We consider the on-line predictive version
  of the standard problem of linear regression;
  the goal is to predict each consecutive response
  given the corresponding explanatory variables
  and all the previous observations.
  We are mainly interested in prediction intervals rather than point predictions.
  The standard treatment of prediction intervals in linear regression analysis
  has two drawbacks:
  (1) the classical prediction intervals
  guarantee that the probability of error
  is equal to the nominal significance level $\epsilon$,
  but this property per se does not imply that the long-run frequency of error
  is close to $\epsilon$;
  (2) it is not suitable for prediction of complex systems
  as it assumes that the number of observations
  exceeds the number of parameters.
  We state a general result showing that in the on-line protocol
  the frequency of error for the classical prediction intervals
  does equal the nominal significance level,
  up to statistical fluctuations.
  We also describe alternative regression models
  in which informative prediction intervals can be found
  before the number of observations exceeds the number of parameters.
  One of these models,
  which only assumes that the observations are independent and identically distributed,
  is popular in machine learning but
  greatly underused in the statistical theory of regression.
\end{abstract}

\section{Introduction}
\label{sec:introduction}

Let $y_n$, $n=1,2,\ldots$, be the sequence of response variables to be predicted,
and let $\mathbf{x}_n=(x_{n,1},\ldots,x_{n,K})$, $n=1,2,\ldots$,
be the corresponding vectors of explanatory variables.
The standard assumption of linear regression analysis
is that the explanatory vectors $\mathbf{x}_n$ are deterministic and
\begin{equation}\label{eq:GL1}
  y_n
  =
  \alpha
  +
  \boldsymbol{\beta}\cdot\mathbf{x}_n
  +
  \xi_n,
\end{equation}
where $\alpha$ is an unknown coefficient,
$\boldsymbol{\beta}\in\bbbr^K$ is an unknown vector of coefficients,
and $\xi_n$, $n=1,2,\ldots$,
are IID (independent and identically distributed) Gaussian random variables
with mean $0$ and unknown variance $\sigma^2>0$
(we will write $\xi_n\sim N(0,\sigma^2)$).
The model (\ref{eq:GL1}) will be called the \emph{Gauss linear model}.
It is the standard textbook model.

The standard classes of problems associated with the Gauss linear model
are parameter estimation,
testing hypotheses about parameters,
and prediction.
In this paper we will be concerned only with prediction,
mainly in the form of prediction intervals rather than point predictions.
(It is natural to concentrate on prediction
as one of the models that we consider, the IID model,
is non-parametric.)

A major drawback of the Gauss linear model
is that the corresponding prediction intervals are uninformative
(i.e., coincide with the whole real line)
unless the number of observations exceeds the number of parameters.
The responses of a complex system
cannot be realistically expected to be modelled
using a small number of parameters,
whereas the number of observations can be very limited.
This motivates consideration of three other models in this paper,
none of which requires that the number of observations
should exceed the number of parameters.

Perhaps the most important of these models
is what we call the \emph{IID model}:
it is only assumed that the sequence of pairs $(\mathbf{x}_n,y_n)$ is IID.
This model is non-parametric,
effectively involving infinitely many parameters.
Despite this,
the model does allow one to obtain informative prediction intervals.
The IID model, however,
also has a fundamental limitation:
informative prediction intervals
become possible only when the number of observations reaches $1/\epsilon$,
where $\epsilon$ is the chosen significance level.

Our third regression model
combines the assumption (\ref{eq:GL1})
with the assumption that $\mathbf{x}_n$
are independent (between themselves and of $\xi_1,\xi_2,\ldots$)
and identically distributed Gaussian random vectors.
We call it the \emph{MVA model},
with MVA referring to ``multivariate analysis''.
It has also been widely discussed
in the statistical literature;
e.g.,
Sampson's \cite{sampson:1974} ``two regressions''
refers to the Gauss linear model and the MVA model.
This model is narrower than both Gauss linear and IID models,
and its strong assumptions ensure that informative prediction intervals
can be produced almost right away.

Finally,
we consider the combination
of the Gauss linear and IID models,
which we call the \emph{IID--Gauss model}:
in addition to (\ref{eq:GL1})
we assume that the explanatory vectors $\mathbf{x}_n$, $n=1,2,\ldots$,
are random and IID
(not necessarily Gaussian, as in the MVA model)
and that the sequence $\xi_1,\xi_2,\ldots$
is independent of the explanatory vectors.
This model, however,
appears to be of secondary importance.
Empirically,
it allows informative prediction intervals
at significance level $\epsilon$
soon after the number of observations exceeds
the minimum of $1/\epsilon$ and the number of parameters.

\begin{figure}[bt]
  \centering
\unitlength 0.50mm
\begin{picture}(150,65)
\put(52,0){\framebox(46,10)[cc]{\textbf{MVA model}}}
\put(27,25){\framebox(96,10)[cc]{IID--Gauss model}}
\put(14,50){\framebox(50,10)[cc]{\textbf{IID model}}}
\put(80,50){\framebox(70,10)[cc]{\textbf{Gauss linear model}}}
\put(75,25){\vector(0,-1){15}}
\put(48,50){\vector(0,-1){15}}
\put(102,50){\vector(0,-1){15}}
\end{picture}
  \caption{\label{fig:models}The four models considered in this paper
    (the three main models are given in boldface).}
\end{figure}
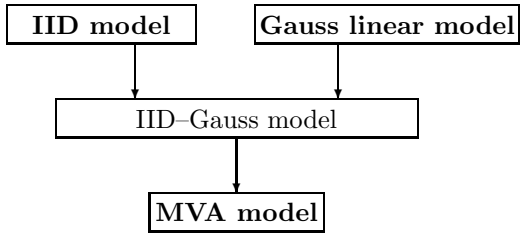

All the models considered in this paper are shown in Figure \ref{fig:models}.
In this paper we begin (in Section \ref{sec:IID})
with the IID model.
This is the most common model used in modern day statistics
and it does not involve the often unrealistic assumption
that the noise variables $\xi_n$ are Gaussian
or that the explanatory vectors $\mathbf{x}_n$ are Gaussian.
An important advantage of the classical Gauss linear model,
considered in Section \ref{sec:GL},
is that the explanatory vectors are not assumed to be IID
(in other words, no ``random design'' is assumed).
This model is essentially equivalent
to making no assumptions whatsoever about the distribution of $\mathbf{x}_n$
and assuming that the $\xi_n$ in (\ref{eq:GL1})
are IID and distributed as $N(0,\sigma^2)$ conditional on
$\mathbf{x}_1,\mathbf{x}_2,\ldots$\,.
The Gauss linear model (understood in this way) and the IID model
are not comparable between themselves,
but both contain the other two models:
the IID--Gauss model (Section \ref{sec:IID-Gauss}),
which is the intersection of the IID and Gauss linear models,
and the MVA model (Section \ref{sec:MVA}),
which makes the further assumption that the explanatory vectors are Gaussian.
These relationships are shown in Figure \ref{fig:models} with the arrows
leading from more general to more specific models.

Fisher (\cite{fisher:1973}, Section IV.3)
emphatically defended the use of the Gauss linear model
even in the case where the distribution of the explanatory vectors is known
(with or without parameters).
There is also a view in the literature
that the Gauss linear model and the MVA model are ``essentially equivalent''
(for a review of some results in this direction,
see \cite{sampson:1974}).
Our conclusion, however, is similar to Brown's \cite{brown:1990}:
when the MVA model is true,
it can be far more useful for prediction;
in particular, it can start giving informative prediction intervals
long before the number of observations reaches
the number of parameters $K$
(or the inverse significance level $1/\epsilon$).

This paper uses a general method of prediction
called conformal prediction.
The method is reviewed in detail in the monograph \cite{vovk/etal:2005book}
and introduced in the work leading up to that monograph.
For each of the four models in Figure \ref{fig:models}
we define a suitable confidence predictor,
i.e., a strategy for producing prediction intervals
or, more generally, prediction regions.
For the IID model we follow \cite{vovk/etal:2005book}
and for the Gauss linear model we use Fisher's classical confidence predictor.
The confidence predictors for the MVA and IID--Gauss models are new.

We are interested in two criteria of quality of confidence predictors,
which we call ``validity'' and ``accuracy''.
For valid confidence predictors,
the probability of error equals the nominal significance level $\epsilon$
(or at least never exceeds $\epsilon$,
in which case we will refer to them as ``conservatively valid'',
or just ``conservative'',
confidence predictors).
(We will not be interested
in approximate versions of validity and conservative validity;
the main problem with such approximate versions
is that it is typically difficult to describe
conditions under which the approximation is good or tolerable.)
The second criterion is applied only to valid confidence predictors:
we want the prediction intervals to be as narrow as possible;
in this paper we, somewhat arbitrarily,
measure the narrowness of a prediction interval $[a,b]$
by its length $b-a$.
In particular,
we want the prediction intervals to become bounded as soon as possible.

Correspondingly, this paper uses two kinds of entities
that one might want to call ``models''.
The first kind is ``hard models'',
such as the four models in Figure \ref{fig:models}.
These are the usual statistical models:
our working hypothesis is that the data set
was generated by one of the probability distributions in the model.
In particular,
the validity of our confidence predictors is allowed to depend on the hard model.
By default, the word ``model'' means ``hard model''.

In addition to the accepted hard model,
one often has other \emph{a priori} information
about the data-generating distribution:
e.g., only a few parameters might provide
the bulk of the information relevant to prediction.
Whereas we might hesitate to include such \emph{a priori} information
in the hard model explicitly,
since it might destroy the validity of our confidence predictor
if this information happened to be far from the truth,
we might still be able to use such information in designing
accurate confidence predictors
provided our model is flexible enough.
A running example in this paper,
introduced in Section \ref{sec:data},
will be a linear system with $100$ parameters
ten of which are felt to be especially important.
This will be our ``soft model''
(not defined formally);
whether it is true or not affects only the accuracy,
but not validity,
of our confidence predictors.

Separation of the available information about the data-generating distribution
into the hard model and soft model
increases robustness of confidence predictors
with respect to modelling errors.
If such an error occurs in the soft model,
the validity of predictions is not affected.
At worst the predictions will become useless,
but they will not become misleading
(with high probability under any distribution in the hard model).
For a further discussion and empirical study,
see \cite{gammerman/vovk:2007local}, Section 4.

The property of validity of conformal predictors 
can be stated in an especially strong form
in the on-line prediction protocol.
It turns out that the true responses fall outside
the corresponding prediction regions independently
for different observations.
In combination with the law of large numbers this implies that,
with high probability,
the frequency of error is approximately equal
to the nominal significance level.
For the classical prediction intervals in the Gauss linear model
this property had been known
(as we will see in Section \ref{sec:GL})
but its importance had not been emphasized
prior to the work leading up to \cite{vovk/etal:2005book}.


For reviews of the theory of conformal prediction see
\cite{gammerman/vovk:2007local} and \cite{shafer/vovk:2008}.
Parts of these two papers are devoted to regression problems.

Section \ref{sec:on-lineI} formally introduces the on-line prediction protocol,
with a more detailed discussion postponed until Section \ref{sec:on-lineII}.
In Section \ref{sec:compression} we describe the method of conformal prediction
and state two key results (proved in Appendix~A):
one asserts the strong validity
and the other universality of conformal predictors.
Section \ref{sec:data} describes an artificial data set
used in later sections for illustrating
the performance of various conformal predictors.
The following four sections,
\ref{sec:IID}--\ref{sec:IID-Gauss},
apply the method of conformal prediction
to the IID, Gauss linear, MVA, and IID--Gauss models,
in this order.
Section \ref{sec:conclusion} concludes.

This version of the paper is an extension of the journal version
\cite{vovk/etal:2009local}.
In addition to the journal version,
it contains a discussion of Fisher's fiducial prediction
(the first subsection of Section \ref{sec:on-lineII})
and three extra appendices, B (containing explicit algorithms),
C (describing an R package implementing those algorithms), and D.

\section{On-line protocol, part I}
\label{sec:on-lineI}

In our prediction protocol,
the task is to sequentially predict $y_n$, $n=1,2,\ldots$,
from $\mathbf{x}_n$ and $(\mathbf{x}_i,y_i)$, $i=1,\ldots,n-1$.
This on-line protocol is popular in machine learning
(see, e.g., \cite{cesabianchi/lugosi:2006} and references therein),
but most statistical research
(except some work on sequential analysis)
is still done in the ``off-line'',
or ``batch'', framework,
where one starts from a complete sample
$(\mathbf{x}_1,y_1),\ldots,(\mathbf{x}_N,y_N)$.
One of the few statisticians advocating the on-line protocol
(under the name ``prequential'', or predictive sequential)
has been Philip Dawid \cite{dawid:1984}.

\subsection*{Weak and strong validity and median accuracy}

To explain what precisely we mean by validity and accuracy,
the two criteria of predictive performance mentioned in Section \ref{sec:introduction},
we will need the notation introduced
in the following description of the on-line prediction protocol.

\bigskip

\noindent
\textsc{On-line prediction protocol}

\medskip

\parshape=7
\IndentI  \WidthI
\IndentII \WidthII
\IndentII \WidthII
\IndentII \WidthII
\IndentII \WidthII
\IndentII \WidthII
\IndentI  \WidthI
\noindent
  FOR $n=1,2,\ldots$:\\
    Predictor observes $\mathbf{x}_n\in\bbbr^K$;\\
    Predictor outputs $\Gamma_n^{\epsilon}\subseteq\bbbr$ for all $\epsilon\in(0,1)$;\\
    Predictor observes $y_n\in\bbbr$;\\
    $\err_n^{\epsilon}:=\III_{y_n\notin\Gamma_n^{\epsilon}}$ for all $\epsilon\in(0,1)$;\\
    $L_n^{\epsilon}:=\sup\Gamma_n^{\epsilon}-\inf\Gamma_n^{\epsilon}$ for all $\epsilon\in(0,1)$\\
  END FOR.

\bigskip

\noindent
(As usual, $\III_F$ is defined to be $1$ if the condition $F$ holds and $0$ if not.)
At each step and for each significance level $\epsilon$,
Predictor outputs a \emph{prediction region}
(usually, although not necessarily, an interval)
$\Gamma_n^{\epsilon}\subseteq\bbbr$.
We require that, for all $n$, the family $\Gamma_n^{\epsilon}$
of prediction regions should be nested:
$\Gamma_n^{\epsilon_1}\subseteq\Gamma_n^{\epsilon_2}$
whenever $\epsilon_1>\epsilon_2$.
An error is registered, $\err_n^{\epsilon}=1$,
if the prediction region fails to contain the true response $y_n$,
and the accuracy of this particular prediction is measured by the length
$L_n^{\epsilon}$
of the corresponding \emph{prediction interval} $\co\Gamma_n^{\epsilon}$
($\co E$ standing for the convex hull of the set $E$).

Let $\Err_n^{\epsilon}:=\err_1^{\epsilon}+\cdots+\err_n^{\epsilon}$
be the cumulative number of errors made up to, and including, step $n$.
In the following sections,
we will find it convenient to distinguish
between two notions of validity,
``weak validity'' and ``strong validity''.

\begin{definition}
  A \emph{confidence predictor} is defined to be
  a measurable prediction strategy
  $
    \Gamma_n^{\epsilon}
    =
    \Gamma^{\epsilon}(\mathbf{x}_1,y_1,\ldots,\mathbf{x}_{n-1},y_{n-1},\mathbf{x}_n)
  $
  in the on-line prediction protocol.
\end{definition}

\begin{definition}
  A confidence predictor is \emph{weakly valid} in some statistical model
  if the probability that $\err_n^{\epsilon}=1$ is $\epsilon$,
  for each $\epsilon\in(0,1)$ and each $n$
  under any probability distribution in the model.
\end{definition}

The definition of weak validity is standard:
cf.\ \cite{cox/hinkley:1974}, (75) on p.~243.
Weak validity by itself does not imply
that $\Err_n/n$ is likely to be close to $\epsilon$
for large $n$.

\begin{definition}
  A confidence predictor is \emph{strongly valid}
  if it is weakly valid and,
  for each $\epsilon\in(0,1)$,
  the events $\err_n^{\epsilon}=1$, $n=1,2,\ldots$,
  are independent.
\end{definition}

Figure \ref{fig:RRCM-errors} below shows the plot of $\Err^{\epsilon}_n$ against $n$
for a specific confidence predictor considered in this paper;
it is typical of our predictors that the slopes of the plots of $\Err^{\epsilon}_n$
are close to the corresponding significance levels $\epsilon$
(we use the significance levels $5\%$, $1\%$ and $0.5\%$ in all our figures,
represented by the corresponding \emph{confidence levels} $1-\epsilon$ in the legends).
This is the only figure in this paper illustrating
the validity of our confidence predictors:
such figures, in view of the mathematical results guaranteeing validity,
tend to be uninformative.

We will measure the accuracy of the predictions
made for the first $n$ observations
by the median $M_n^{\epsilon}$ of the sequence
$L_1^{\epsilon},\ldots,L_n^{\epsilon}$;
again,
this measure is arbitrary, to a large degree.
A plot of $M_n^{\epsilon}$ against $n$ will be called the \emph{median-accuracy plot};
examples of such plots are given
in Figures \ref{fig:RRCM-med} and \ref{fig:SM-med}--\ref{fig:comb-med}.

Unfortunately,
the simple notions of validity introduced earlier
have to be extended to become useful for our purpose.
This is needed because, e.g.,
the classical prediction intervals are uninformative
before the number of observations reaches the number of parameters,
and so for small $n$ the error probability is zero rather than $\epsilon$.
Let $\NNN$ be a set of positive integer numbers
(we are mainly interested in the case where $\NNN$ has the form $\{m,m+1,\ldots\}$).

\begin{definition}
  We say that a confidence predictor is \emph{weakly valid for} $n\in\NNN$
  in a statistical model if the probability is $\epsilon$
  that it makes an error, $\err_n^{\epsilon}=1$, at step $n$
  under any probability distribution in the model
  and for all $n\in\NNN$ and $\epsilon\in(0,1)$.
  It is \emph{strongly valid for} $n\in\NNN$ if, in addition,
  $\err_n^{\epsilon}$, $n\in\NNN$,
  are independent for any fixed $\epsilon$.
\end{definition}

\subsection*{The role of the on-line protocol}

The exposition of this paper is based on the on-line protocol,
but the majority of our findings are not constrained to this specific protocol.
For example,
the fact that valid and informative prediction intervals can become feasible
in the MVA model
before the number of observations exceeds the number of parameters
does not depend on the prediction protocol.
In the absence of the on-line protocol, however,
``validity'' should be understood in the standard sense of weak validity.

\section{Conformal prediction}
\label{sec:compression}


In this section we
define a class of confidence predictors,
called conformal predictors,
and state results about their validity and universality,
in a certain sense.

\subsection*{Notions of sufficiency}

Fix some \emph{observation space} $\Zeta$.
We will be interested in the space $\Zeta=\bbbr^{K}\times\bbbr$
of pairs $(\mathbf{x},y)$;
in general, $\Zeta$ is a measurable space
assumed to be Luzin,
to ensure the existence of regular conditional probabilities.
To define conformal predictors,
we will need not only a statistical model on $\Zeta^{\infty}$
but also a sequence of sufficient statistics $S_n:\Zeta^n\to\Sigma_n$,
$n=1,2,\ldots$;
we will always assume that $\Sigma_n=S_n(\Zeta^n)$.
We will need a strengthened form of sufficiency;
in our definitions we mainly follow Lauritzen
\cite{lauritzen:1988}, Section II.2.

The sequence $(S_n)$ is \emph{algebraically transitive}
if there exists a sequence of measurable functions
$F_n:\Sigma_{n-1}\times\Zeta\to\Sigma_n$,
$n=2,3,\ldots$,
such that
\begin{equation*}
  S_n(\zeta_1,\ldots,\zeta_{n-1},\zeta_n)
  =
  F_n(S_{n-1}(\zeta_1,\ldots,\zeta_{n-1}),\zeta_n)
\end{equation*}
for all $(\zeta_1,\ldots,\zeta_{n-1},\zeta_n)\in\Zeta^n$.
Intuitively,
$S_n(\zeta_1,\ldots,\zeta_n)$ is the summary of the first $n$ observations,
and the condition of algebraic transitivity
means that the summary can be updated on-line.

The sequence $(S_n)$ is \emph{totally sufficient} 
for a statistical model $\PPP$ on $\Zeta^{\infty}$
if, for each $n=1,2,\ldots$:
\begin{itemize}
\item
  $S_n$ is sufficient for $\PPP$;
\item
  $\zeta_1,\ldots,\zeta_n$ and $\zeta_{n+1},\zeta_{n+2},\ldots$
  are conditionally independent given $S_n(\zeta_1,\ldots,\zeta_n)$,
  where $(\zeta_1,\zeta_2,\ldots)\sim P$, for any $P\in\PPP$.
\end{itemize}
The second condition ensures that $S_n(\zeta_1,\ldots,\zeta_n)$
carries all information in $\zeta_1,\ldots,\zeta_n$
that can be used for predicting the future observations $\zeta_{n+1},\zeta_{n+2},\ldots$\,.

A sequence of statistics that is both algebraically transitive
and totally sufficient will be called an \emph{ATTS sequence}.
In the rest of this paper we will often say ``model''
to mean a statistical model $\PPP$ equipped with an ATTS sequence $(S_n)$.
This makes the word ``model'' ambiguous
as we often omit ``statistical'' in ``statistical model'',
but this should not lead to misunderstandings.


Each of the four statistical models considered in this paper
(see Figure \ref{fig:models})
will be complemented with an ATTS sequence;
in all four cases
the observation space $\Zeta$ will be $\bbbr^{K}\times\bbbr$.

\subsection*{Testing conformity}

The main ingredient of conformal prediction
is statistical testing of conformity of a new observation $\zeta_n$
to the old observations $\zeta_1,\ldots,\zeta_{n-1}$.
In general, our statistical tests will be randomized.

Fix a statistical model $\PPP$ with an ATTS sequence $S_n:\Zeta^n\to\Sigma_n$.
Define $\Sigma_0$ to be a fixed one-element set.
Any sequence of measurable functions
$A_n:\Sigma_{n-1}\times\Zeta\to\bbbr$,
$n=1,2,\ldots\,$,
is called a \emph{nonconformity measure};
$A_n$ will be our test statistics.
Given a nonconformity measure $(A_n)$,
for each sequence $\zeta_1,\zeta_2,\ldots$ of observations
and each sequence $\tau_1,\tau_2,\ldots\in[0,1]^{\infty}$
we define the \emph{$p$-values}
\begin{multline}\label{eq:p}
  p_n
  =
  p_n(\zeta_1,\ldots,\zeta_n,\tau_n)
  :=\\
  \Prob
  \left(
    A_n\rnd > A_n\obs
    \given
    S_n\rnd = S_n\obs
  \right)
  +
  \tau_n
  \Prob
  \left(
    A_n\rnd = A_n\obs
    \given
    S_n\rnd = S_n\obs
  \right),
  \\
  n=1,2,\ldots,
\end{multline}
where
$A_n\rnd:=A_n(S_{n-1}(\xi_1,\ldots,\xi_{n-1}),\xi_n)$
and
$S_n\rnd:=S_n(\xi_1,\ldots,\xi_n)$
are the ``random'' values,
$A_n\obs:=A_n(S_{n-1}(\zeta_1,\ldots,\zeta_{n-1}),\zeta_n)$
and
$S_n\obs:=S_n(\zeta_1,\ldots,\zeta_n)$
are the ``observed'' values,
and the probabilities are taken with respect to
$(\xi_1,\xi_2,\ldots)\sim P$ for some $P\in\PPP$.
Since $S_n$ are sufficient statistics,
$p_n$ do not depend on $P\in\PPP$
(at least for a suitable choice of regular conditional probabilities).
We will be interested in two cases:
\emph{deterministic}, where $\tau_n=1$ for all $n$,
and \emph{randomized},
where $\tau_1,\tau_2,\ldots$ are generated independently
from the uniform distribution $U$ on $[0,1]$
(such $\tau_1,\tau_2,\ldots$ model the output of a random numbers generator).
\begin{theorem}\label{thm:1}
  Suppose that the sequence of observations $(\zeta_1,\zeta_2,\ldots)\in\Zeta^{\infty}$
  is generated from a probability distribution $P\in\PPP$
  and that the random numbers $(\tau_1,\tau_2,\ldots)\sim U^{\infty}$
  are independent of the observations.
  The $p$-values (\ref{eq:p}) are then independent
  and distributed uniformly on $[0,1]$:
  \begin{equation*}
    (p_1,p_2,\ldots) \sim U^{\infty}.
  \end{equation*}
\end{theorem}
For a proof of this theorem,
see Appendix~A.
The fact that $p_n\sim U$ is well known,
at least in the continuous case
(see, e.g., \cite{cox/hinkley:1974}, p.~66;
(\ref{eq:p}) is a version of (1) in \cite{cox/hinkley:1974}).

\subsection*{Conformal prediction}

We start by extending,
and spelling out in a greater detail,
the notion of a confidence predictor:
in the general theory of this section
and in its application to the IID model in Section \ref{sec:IID}
we will need an element (typically quite small) of randomization
in confidence predictors.

\begin{definition}
  A \emph{randomized confidence predictor}
  is a measurable function
  which maps every significance level $\epsilon\in(0,1)$,
  every data sequence $\mathbf{x}_1,y_1,\ldots,\mathbf{x}_{n-1},y_{n-1}$,
  every vector $\mathbf{x}_n$ of explanatory variables,
  and every number $\tau\in[0,1]$
  to a set
  $
    \Gamma_n^{\epsilon}
    =
    \Gamma^{\epsilon}(\mathbf{x}_1,y_1,\ldots,\mathbf{x}_{n-1},y_{n-1},\mathbf{x}_n,\tau)
    \subseteq
    \bbbr
  $.
  We will use the notation $\Gamma_n^{\epsilon}$ when the data sequence,
  the vector of explanatory variables, and the number $\tau$ are clear from the context.
\end{definition}

Let the observation space be $\Zeta=\bbbr^{K}\times\bbbr$.
Once the $p$-values (\ref{eq:p}) are defined,
we can use them for confidence prediction
(this is a standard procedure;
cf.\ \cite{cox/hinkley:1974}, (76) on p.~243):
we set
\begin{multline}\label{eq:Gamma}
  \Gamma^{\epsilon}(\mathbf{x}_1,y_1,\ldots,\mathbf{x}_{n-1},y_{n-1},\mathbf{x}_n,\tau_n)\\
  {}:=
  \left\{
    y\in\bbbr
    \st
    p_n((\mathbf{x}_1,y_1),\ldots,(\mathbf{x}_{n-1},y_{n-1}),(\mathbf{x}_n,y),\tau_n)
    >
    \epsilon
  \right\}.
\end{multline}

\begin{definition}
  The randomized confidence predictor
  defined by (\ref{eq:Gamma})
  is called
  the \emph{smoothed conformal predictor determined by}
  the nonconformity measure $(A_n)$.
  A \emph{smoothed conformal predictor} is
  a smoothed conformal predictor determined by some nonconformity measure.
\end{definition}

The following statement immediately follows from Theorem \ref{thm:1}
and asserts that smoothed conformal predictors are strongly valid.

\begin{corollary}\label{cor:validity}
  If the sequence of observations $(\mathbf{x}_n,y_n)$, $n=1,2,\ldots$,
  is generated by a probability distribution $P\in\PPP$
  and a smoothed conformal predictor
  is fed with random numbers $(\tau_1,\tau_2,\ldots)\sim U^{\infty}$
  independent of the observations,
  the error sequence $\err_1^{\epsilon},\err_2^{\epsilon},\ldots$
  at any significance level $\epsilon$
  is a sequence of IID Bernoulli random variables with parameter $\epsilon$.
\end{corollary}

The adjective ``smoothed'' refers to using random numbers;
if we take $\tau_n=1$ for all $n=1,2,\ldots\,$,
we will obtain the definition of a ``deterministic conformal predictor'',
or just ``conformal predictor'',
and in this case we omit $\tau_n$ from our notation.

\begin{definition}
  A \emph{conformal predictor} is the confidence predictor
  defined by
  \begin{multline*}
    \Gamma^{\epsilon}(\mathbf{x}_1,y_1,\ldots,\mathbf{x}_{n-1},y_{n-1},\mathbf{x}_n)\\
    {}:=
    \left\{
      y\in\bbbr
      \st
      p_n((\mathbf{x}_1,y_1),\ldots,(\mathbf{x}_{n-1},y_{n-1}),(\mathbf{x}_n,y),1)
      >
      \epsilon
    \right\},
  \end{multline*}
  where the $p$-values $p_n$ are defined by (\ref{eq:p}).
\end{definition}

Notice that when a conformal predictor makes an error,
the corresponding smoothed conformal predictor also makes an error.
In combination with Corollary~\ref{cor:validity},
we can see that conformal predictors are \emph{conservative},
in the sense that, for each $\epsilon$,
their error sequence $\err_1^{\epsilon},\err_2^{\epsilon},\ldots$
is dominated by a sequence of IID Bernoulli random variables with parameter $\epsilon$.
In particular,
whereas we have
$\lim_{n\to\infty}(\Err^{\epsilon}_n/n)=\epsilon$ a.s.\ for
smoothed conformal predictors,
we only have
$\limsup_{n\to\infty}(\Err^{\epsilon}_n/n)\le\epsilon$ a.s.\ for
conformal predictors.

We will see that there is no difference between conformal predictors
and the corresponding smoothed conformal predictors
for the Gauss linear model and $n\ge K+3$
since the second addend on the right-hand side of (\ref{eq:p})
is then zero.
There is also no difference for the MVA model and $n\ge3$;
however,
the difference is important
(although usually barely noticeable on error and accuracy plots)
for the IID model.


A natural question is whether there are other ways
to achieve validity,
except conformal prediction.
The following theorem will give a negative answer to a version of this question.


\begin{definition}
  A confidence predictor $\Gamma$ is \emph{invariant}
  if $\Gamma_n^{\epsilon}$, $n>1$, depends on the first $n-1$ observations
  only through the value of $S_{n-1}$ on those observations.
\end{definition}

The use of invariant confidence predictors
is natural in view of the sufficiency principle;
see, e.g., \cite{cox/hinkley:1974}, Section 2.3 (ii).
Let $\NNN$ be a set of positive integers.
We say that a confidence predictor $\Gamma^{\dagger}$
is \emph{at least as accurate as} another confidence predictor $\Gamma$
for $n\in\NNN$
if
\begin{equation*}
  (\Gamma^{\dagger})^{\epsilon}(\mathbf{x}_1,y_1,\ldots,\mathbf{x}_{n-1},y_{n-1},\mathbf{x}_n)
  \subseteq
  \Gamma^{\epsilon}(\mathbf{x}_1,y_1,\ldots,\mathbf{x}_{n-1},y_{n-1},\mathbf{x}_n)
\end{equation*}
for all $\epsilon$,
all $n\in\NNN$,
and $P$-almost all $\mathbf{x}_1,y_1,\ldots,\mathbf{x}_{n-1},y_{n-1},\mathbf{x}_n$,
under any probability distribution $P\in\PPP$.

Recall that a statistic $S$ taking values in a measurable space $\Sigma$
is said to be \emph{boundedly complete}
(with respect to the statistical model $\PPP$)
if, for any bounded measurable function $f:\Sigma\to\bbbr$,
the following condition is satisfied:
the expected value $\Expect_P(f(S))$ of $f(S)$
is zero under all $P\in\PPP$
only if $f(S)=0$ $P$-almost surely for all $P\in\PPP$.
\begin{theorem}\label{thm:2}
  Let $\NNN$ be a set of positive integers.
  Suppose the ATTS statistics $S_n$ are boundedly complete for $n\in\NNN$.
  If a confidence predictor $\Gamma$ is invariant and weakly valid for $n\in\NNN$,
  then there is a conformal predictor that is at least as accurate as $\Gamma$
  for $n\in\NNN$.
\end{theorem}
This theorem is also proved in Appendix~A.
An important step towards its proof
was made by Takeuchi (\cite{\Takeuchi}, p.~31).

The condition of bounded completeness holds
for the Gauss linear model and the MVA model
by the standard completeness result for exponential statistical models
(see, e.g., Theorem~4.1 in \cite{lehmann:1986}),
and it is also known to hold for the IID model
(see the theorem on p.~797 in \cite{bell/etal:1960}
or Theorem~1 in \cite{mattner:1996}).

\section{Data set}
\label{sec:data}

We will illustrate the accuracy of various confidence predictors
using the following artificially generated data set
with $600$ observations and $K=100$ explanatory variables.
The components $x_{n,k}$ of $\mathbf{x}_n$
are independently generated from $N(0,1)$,
and the responses $y_n$ are generated according to (\ref{eq:GL1})
with $\xi_n\sim N(0,1)$
independent between themselves and of all $x_{n,k}$,
with $\alpha=100$
and with the following components $\beta_k$ of $\boldsymbol{\beta}$:
\begin{equation*}
  \beta_k
  :=
  \begin{cases}
    (-1)^{k-1} 10 & k=1,\ldots,10\\
    (-1)^{k-1} & k=11,\ldots,100.
  \end{cases}
\end{equation*}

The probability distribution generating this data set
belongs to all four models considered in this paper
(Figure \ref{fig:models}).
It is natural to expect that more specific models,
when true,
will lead to better predictions.
In one respect this is true:
more general models allow informative predictions later,
as shown in Table \ref{tab:limitations}
(to be explained in later sections).
However, soon after the threshold given in the table is reached,
the quality of prediction becomes very similar on our data set.

\begin{table}
\caption{Steps at which informative prediction becomes possible for the four models;
  $\epsilon$ is the significance level ($\epsilon<1/2$ is assumed)
  and $K$ is the number of parameters.}
\label{tab:limitations}
\begin{tabular}{cc}
\hline
& The first step at which prediction intervals\\Model & can become informative\\
\hline
IID model & $\lceil 1/\epsilon\rceil$\\
Gauss linear model & $K+3$\\
MVA model & 3\\
IID--Gauss model & $\min\left(\lceil 1/\epsilon\rceil,K+3\right)$\\
\hline
\end{tabular}
\end{table}

The (informal) soft model guiding the choice of the nonconformity measure
will include the assumption of linearity (\ref{eq:GL1})
and the knowledge, or guess,
that the first 10 explanatory variables
are much more important than the rest.

Relationship (\ref{eq:GL1}) between the response and explanatory variables
can be written as
\begin{equation}\label{eq:GL2}
  y_n
  =
  \boldsymbol{\gamma}\cdot\mathbf{z}_n
  +
  \xi_n,
\end{equation}
where
\begin{equation*}
  \boldsymbol{\gamma}
  :=
  \begin{pmatrix}
    \alpha\\
    \boldsymbol{\beta}
  \end{pmatrix}
  \in
  \bbbr^{K+1}
  \text{ and }
  \mathbf{z}_n
  :=
  \begin{pmatrix}
    1\\
    \mathbf{x}_n
  \end{pmatrix}
  \in
  \bbbr^{K+1}.
\end{equation*}
For $l=1,2,\ldots$,
let $\mathbfit{Z}_l$ be the $l\times(K+1)$ matrix
whose rows are $\mathbf{z}'_i$, $i=1,\ldots,l$,
and $\mathbf{y}_l$ be the vector whose $i$th element is $y_i$, $i=1,\ldots,l$.
We will sometimes refer to the first column of $\mathbfit{Z}_l$
as the \emph{dummy} column.

\section{The IID model}
\label{sec:IID}

The statistical model considered in this section is non-parametric:
we simply assume that the observations $(\mathbf{x}_n,y_n)$ are IID.
Notice that this does not involve the assumption of linearity
of the ``true'' regression function
or the assumption of a Gaussian noise.
Linearity is, however, an important component of the soft model
used for choosing a suitable nonconformity measure.

The ATTS statistics are
\begin{equation*}
  S_n
  :=
  \lbag
    (\mathbf{x}_1,y_1),
    \ldots,
    (\mathbf{x}_n,y_n)
  \rbag,
\end{equation*}
where we use $\lbag a_1,\ldots,a_n\rbag$ to denote the bag,
or multiset, consisting of $a_1,\ldots,a_n$
(some of these elements may coincide).
For each $n$,
the conditional distribution of $(\xi_1,\ldots,\xi_n)$ given that
\begin{equation*}
  \lbag
    \xi_1,
    \ldots,
    \xi_n
  \rbag
  =
  \lbag
    (\mathbf{x}_1,y_1),
    \ldots,
    (\mathbf{x}_n,y_n)
  \rbag,
\end{equation*}
where $\xi_i$ are IID random elements taking values in $\bbbr^K\times\bbbr$,
coincides (with probability one) 
with the probability distribution on the orderings
$(\mathbf{x}_{\pi(1)},y_{\pi(1)}),\ldots,(\mathbf{x}_{\pi(n)},y_{\pi(n)})$
of the bag $\lbag(\mathbf{x}_1,y_1),\ldots,(\mathbf{x}_n,y_n)\rbag$
that arises when different permutations $\pi$ are chosen
with the same probability $1/n!$.

The IID model is typical in that
there is a great flexibility in choosing a nonconformity measure
for use in conformal prediction.
Suppose, e.g., that the number of explanatory variables $K$
is too large for us to estimate all the $\beta_k$ and $\alpha$
in the soft model (\ref{eq:GL1}).
We believe, however, that the first $K^{\dagger}_n \ll K$ of the explanatory variables
are especially important,
and it is feasible to estimate
the corresponding $\beta_k$, $k=1,\ldots,K^{\dagger}_n$, and $\alpha$.

Fix temporarily a positive integer number $n$.
We will write $\mathbf{y}$ for $\mathbf{y}_n$,
$\mathbfit{Z}$ for $\mathbfit{Z}_n$
and $K^{\dagger}$ for $K^{\dagger}_n$.
Let $\mathbfit{U}$ be the submatrix of $\mathbfit{Z}$
consisting of the first $K^{\dagger}+1$ columns of $\mathbfit{Z}$:
those that correspond to the explanatory variables deemed to be useful at this stage
plus the dummy column $\boldsymbol{1}$.
To test the conformity of the $n$th observation
to the first $n-1$ observations,
we will first fit a hyperplane
to all $n$ observations using the relevant explanatory variables.
Applying a small ``ridge coefficient'' $a>0$
to avoid the need to invert singular matrices,
we obtain the vector of residuals
\begin{equation}\label{eq:residuals1}
  \mathbf{e}
  :=
  \mathbf{y}
  -
  \mathbfit{U}
  \left(
    \mathbfit{U}'\mathbfit{U}+a\mathbfit{I}
  \right)^{-1}
  \mathbfit{U}'\mathbf{y},
\end{equation}
whose components will be denoted $e_1,\ldots,e_n$.

We will be interested in the conformal predictor
determined by the nonconformity measure
\begin{equation}\label{eq:nonconformity-de-Finetti}
  A_n
  \left(
    S_{n-1}
    \left(
      \mathbf{x}_1,y_1,
      \ldots,
      \mathbf{x}_{n-1},y_{n-1}
    \right),
    (\mathbf{x}_{n},y_{n})
  \right)
  :=
  \lvert e_n\rvert.
\end{equation}
Deleted and, especially, studentized residuals
would also be a natural choice
(see, e.g., \cite{vovk/etal:2005book}, pp.~34--35).
In our experience, however, the difference is not significant,
and we stick to the simplest choice.
The confidence predictor obtained from this conformal predictor
by replacing the prediction regions $\Gamma_n^{\epsilon}$
with the prediction intervals $\co\Gamma_n^{\epsilon}$
will be called the \emph{IID predictor}
(cf.\ the comments at the end of this section).

The IID predictor can be implemented fairly efficiently.
First notice that for the IID model
the formula (\ref{eq:p}) for $p$-values
can be simplified to
\begin{equation}\label{eq:p-de-Finetti}
  p_n
  =
  \frac
  {
    \left|
      \left\{
        i \st \alpha_i > \alpha_n
      \right\}
    \right|
    +
    \tau_n
    \left|
      \left\{
        i \st \alpha_i = \alpha_n
      \right\}
    \right|
  }
  {
    n
  },
\end{equation}
where $\alpha_i:=A_n(\lbag\zeta_1,\ldots,\zeta_{i-1},\zeta_{i+1},\ldots,\zeta_n\rbag,\zeta_i)$,
$i$ ranges over $\{1,\ldots,n\}$,
and $\lvert E\rvert$ stands for the size of the set $E$.
In the case of the nonconformity measure
(\ref{eq:nonconformity-de-Finetti}),
$\alpha_i=\lvert e_i\rvert$.
The residuals (\ref{eq:residuals1}) can be written in the form
\begin{equation*}
  \mathbf{e}
  =
  \mathbf{y}
  -
  \mathbfit{U}
  \left(
    \mathbfit{U}'\mathbfit{U}+a\mathbfit{I}
  \right)^{-1}
  \mathbfit{U}'\mathbf{y}
  =
  \mathbfit{C}\mathbf{y},
\end{equation*}
where $\mathbfit{C}$ is the matrix
$
  \mathbfit{I}
  -
  \mathbfit{U}
  \left(
    \mathbfit{U}'\mathbfit{U}+a\mathbfit{I}
  \right)^{-1}
  \mathbfit{U}'
$,
not depending on the response variables.
If we fix the first $n-1$ response variables $y_i$ and vary the last one, $y$,
the residuals $e_i=e_i(y)$, $i=1,\ldots,n$, become linear functions of $y$
(this fact will also be used in Section \ref{sec:MVA}).
By (\ref{eq:p-de-Finetti}) with $\tau_n:=1$,
the $p$-value is the fraction of $i=1,\ldots,n$
satisfying $\lvert e_i(y)\rvert \ge \lvert e_n(y)\rvert$;
therefore, as $y$ varies from $-\infty$ to $\infty$,
the $p$-value can change only at the at most $2n-2$ points
(called \emph{critical} points)
which are solutions to the linear equations
$e_i(y) = e_n(y)$ and $e_i(y) = -e_n(y)$.
This divides the real line into at most $4n-3$ intervals:
the critical points,
considered as degenerate closed intervals,
the open intervals bounded on both sides by adjacent critical points,
and the two unbounded open intervals to the left of the leftmost critical point
and to the right of the rightmost critical point;
if there are no critical points, this collapses into one unbounded open interval $\bbbr$.
We can compute the $p$-value for one point in each of these intervals
and then compute $\Gamma_n^{\epsilon}$ as the union of the intervals
with $p$-values exceeding $\epsilon$.
The computation of the IID prediction interval
$\co\Gamma_n^{\epsilon}$
can be simplified if we notice that the set $\Gamma_n^{\epsilon}$ is closed
(which is opposite to what we will have for the Gauss linear and MVA models):
assuming that the set of critical points is non-empty,
$\co\Gamma_n^{\epsilon}$ is bounded if and only if the two unbounded intervals
have $p$-values at most $\epsilon$,
in which case the end-points of $\co\Gamma_n^{\epsilon}$
can be found as the leftmost and rightmost critical points
with $p$-values exceeding $\epsilon$.
Computing $\Gamma_n^{\epsilon}$ and $\co\Gamma_n^{\epsilon}$ from scratch
(e.g., without using the results of computations from the previous steps
of the on-line protocol)
takes time $O(n\log n)$
(see \cite{vovk/etal:2005book}, p.~33).

For use in our experiments with the artificial data set
described in Section \ref{sec:data},
we take
\begin{equation}\label{eq:K}
  K^{\dagger}_n
  :=
  \begin{cases}
    10 & \text{if $n<103$}\\
    100 & \text{otherwise},
  \end{cases}
\end{equation}
and so define $\mathbfit{U}$ as the first $11$ columns of $\mathbfit{Z}$ if $n<103$
and as the full $\mathbfit{Z}$ otherwise.
Our chosen value for the threshold, $103$, appeared to us slightly less arbitrary
than other choices,
since it is the first step when the classical prediction intervals
(see (\ref{eq:GammaGLM}) below)
become bounded.
However, the quality of the estimates of $\alpha$
and the $100$ components of $\boldsymbol{\beta}$
is still poor when $n$ is close to $103$.
This affects the quality of our prediction intervals
but does not show on the median-accuracy plots.
The value of the ridge coefficient is always $a=0.01$.

\begin{figure}[bt]
  \centering
    \ifCOLOR
      \makebox{\includegraphics[width=\picturewidth]{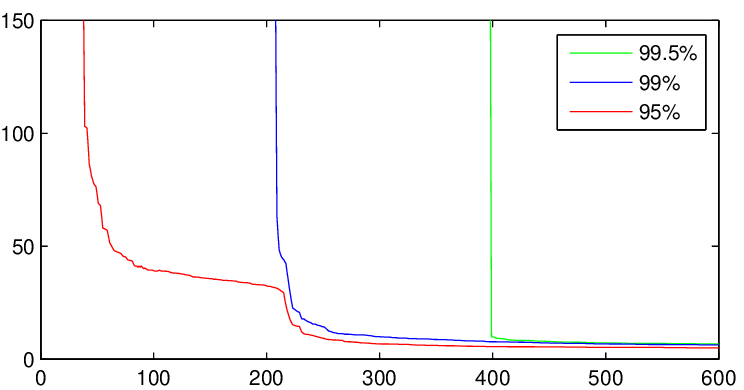}}
    \fi
    \ifnotCOLOR
      \makebox{\includegraphics[width=\picturewidth]{RRCM0_Med_bw.eps}}
    \fi
  \caption{\label{fig:RRCM-med}The median-accuracy plot for the IID predictor.
    The three significance levels used in this and all the following figures
    are $\epsilon=0.05,0.01,0.005$,
    shown in the form $100(1-\epsilon)\%$
    (the corresponding confidence levels) in the legends.}
\end{figure}

As Figure \ref{fig:RRCM-med} shows,
the IID predictor works well for our data set
if the significance level is not too demanding:
it can be seen from (\ref{eq:p-de-Finetti})
(with $\tau_n:=1$)
that for the IID prediction interval $\co\Gamma_n^{\epsilon}$ to be bounded
the number of observations $n$ has to be at least $1/\epsilon$
(as Table \ref{tab:limitations} says).
For example,
for the significance level $\epsilon=0.5\%$,
the IID predictor requires $200$ observations to produce bounded predictions,
and this shows on the median-accuracy plot at $n=399$
(since for $n<399$ at least half of the observed prediction intervals
are infinitely wide).


The IID model is non-parametric
but we can see that it still admits valid confidence predictors
(or conservative confidence predictors
if one insists on using deterministic predictors).
The threshold $1/\epsilon$ can be said to play the role
of the number of parameters,
and the non-parametric nature of the model is reflected in the fact that
$1/\epsilon\to\infty$ as $\epsilon\to0$.
Since $1/\epsilon$ tends to $\infty$ relatively slowly,
such an infinite-dimensional model may be better
for the purpose of prediction
than a $K$-dimensional model with a very large $K$.

Theorem~\ref{thm:2} is not directly applicable to the IID model,
since only smoothed conformal predictors are valid,
as the latter term is used in this paper.
Two results of the same nature about the IID model are stated
in \cite{vovk/etal:2005book}, Section 2.4.

\begin{figure}[bt]
  \centering
    \ifCOLOR
      \makebox{\includegraphics[width=\picturewidth]{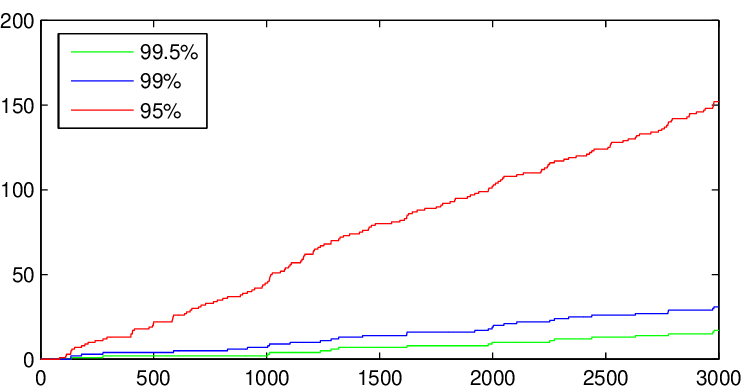}}
    \fi
    \ifnotCOLOR
      \makebox{\includegraphics[width=\picturewidth]{RRCM0_Err_bw.eps}}
    \fi
  \caption{\label{fig:RRCM-errors}The cumulative numbers of errors
    made by the IID predictor:
    $\Err_n^{\epsilon}$ is plotted against $n$.}
\end{figure}

There are two sources of conservativeness
for the IID predictor as described above
(and used for producing Figure \ref{fig:RRCM-med}).
First,
we used a deterministic predictor (taking $\tau_n=1$ for all $n$),
and second,
we replaced each prediction region by its convex hull.
Our experiments (see, e.g., Figure \ref{fig:RRCM-errors})
show that we still have approximate validity.

For each model considered in this paper
except the Gauss linear model
we define a nonconformity measure
involving the matrix $\mathbfit{U}$ defined earlier in this section.
In the case of the IID model,
we have used the nonconformity measure (\ref{eq:nonconformity-de-Finetti})
and called the corresponding conformal predictor
with $\Gamma_n^{\epsilon}$ replaced by $\co\Gamma_n^{\epsilon}$
the IID predictor
(it was called ``Ridge Regression Confidence Machine''
in \cite{vovk/etal:2005book}%
\Extra{, as well as in the original paper where it was introduced}).
Of course,
our brief term is somewhat misleading:
it should always be borne in mind
that the conformal predictor leading to the IID predictor
is only one of many conformal predictors that can be defined in the IID model.
Similarly, in the following three sections
we will introduce the Gauss predictor, the MVA predictor,
and the IID--Gauss predictor,
which will also correspond to specific nonconformity measures.

\section{The Gauss linear model}
\label{sec:GL}

Let
$
  \hat{\boldsymbol{\gamma}}_l
  :=
  (\mathbfit{Z}'_l\mathbfit{Z}_l)^{-1}\mathbfit{Z}'_l\mathbf{y}_l
$
be the least-squares estimate of the parameter vector
$\boldsymbol{\gamma}$ in (\ref{eq:GL2})
from the first $l$ observations.
For simplicity,
we will assume that the matrix $\mathbfit{Z}_{l}$ has full rank
(i.e., $\rank\mathbfit{Z}_l=\min(l,K+1)$)
for all $l$;
this implies that $\hat{\boldsymbol{\gamma}}_l$ is well defined for $l\ge K+1$.

Let $\hat{y}_n$ be the least-squares prediction
$\hat{\boldsymbol{\gamma}}_{n-1}\cdot\mathbf{z}_n$
for $y_n$
and
\begin{equation*}
  \hat{\sigma}^2_{l}
  :=
  \frac{1}{l-K-1}
  (\mathbf{y}_l-\mathbfit{Z}_l\hat{\boldsymbol{\gamma}}_l)'
  (\mathbf{y}_l-\mathbfit{Z}_l\hat{\boldsymbol{\gamma}}_l)
\end{equation*}
be the standard estimate of $\sigma^2$
from $\mathbfit{Z}_{l}$ and $\mathbf{y}_{l}$.
It is well known that
in the Gauss linear model
the ratio
\begin{equation}\label{eq:pivotGLM}
  T_n
  :=
  \frac
  {
    y_n-\hat{y}_n
  }
  {
    \sqrt{1+\mathbf{z}'_n(\mathbfit{Z}'_{n-1}\mathbfit{Z}_{n-1})^{-1}\mathbf{z}_n}
    \hat{\sigma}_{n-1}
  },
  \quad
  n=K+3,K+4,\ldots,
\end{equation}
has the $t$-distribution with $n-K-2$ degrees of freedom.
This gives the classical weakly valid prediction interval for the $n$th response,
\begin{multline}\label{eq:GammaGLM}
  \Gamma_n^{\epsilon}
  :=
    \left\{
      y\in\bbbr
      \st
      \left|
        y-\hat{y}_n
      \right|
      <
      t_{n-K-2}^{\epsilon/2}
      \sqrt{1+\mathbf{z}'_n(\mathbfit{Z}'_{n-1}\mathbfit{Z}_{n-1})^{-1}\mathbf{z}_n}
      \hat{\sigma}_{n-1}
    \right\},\\
    n\ge K+3,
\end{multline}
where $t_{m}^{\delta}$ is the upper $\delta$ point of the $t$-distribution
with $m$ degrees of freedom.
(See, e.g., \cite{seber/lee:2003}, (5.27).)
We set $\Gamma_n^{\epsilon}$ to $\bbbr$ when $n<K+3$.

Later in this section we will see that Corollary~\ref{cor:validity}
implies the following property of the classical prediction intervals
for the Gauss linear model.

\begin{corollary}\label{prop:GLM}
  Let $\epsilon\in(0,1)$.
  The events $y_n\notin\Gamma^{\epsilon}_n$, $n=K+3,K+4,\ldots$,
  are independent.
  In particular,
  the confidence predictor (\ref{eq:GammaGLM}) is strongly valid
  for $n\ge K+3$.
\end{corollary}


\begin{remark}
  Corollary \ref{prop:GLM} and, more generally,
  the fact that the statistics (\ref{eq:pivotGLM}) are independent
  was established in \cite{diaz/etal:1983}
  following \cite{oreilly/quesenberry:1973}.
  It is interesting that both papers use the independence of (\ref{eq:pivotGLM})
  for testing rather than for prediction.
\end{remark}

Let us now see that some conformal predictor
outputs the classical prediction intervals (\ref{eq:GammaGLM}).
This will demonstrate that Corollary~\ref{prop:GLM}
is indeed a special case of Corollary~\ref{cor:validity}.

The ATTS statistics for the Gauss linear model are
\begin{equation*} 
  S_n
  (\mathbf{x}_1,y_1,\ldots,\mathbf{x}_n,y_n)
  :=
  \left(
    \mathbf{x}_1,\ldots,\mathbf{x}_n,
    \sum_{i=1}^n y_i,
    \sum_{i=1}^n y_i\mathbf{x}_i,
    \sum_{i=1}^n y^2_i
  \right).
\end{equation*}
(It is natural to have $\mathbf{x}_1,\ldots,\mathbf{x}_n$ as components of $S_n$,
although they are superfluous under our original definition,
in which $\mathbf{x}_1,\mathbf{x}_2,\ldots$ are deterministic.)
The prediction intervals (\ref{eq:GammaGLM})
are precisely the prediction regions
output by the conformal predictor corresponding to the nonconformity measure
\begin{multline}\label{eq:nonconformityGLM}
  A_n
  \left(
    S_{n-1}
    \left(
      \mathbf{x}_1,y_1,
      \ldots,
      \mathbf{x}_{n-1},y_{n-1}
    \right),
    (\mathbf{x}_{n},y_{n})
  \right)\\
  :=
  \frac
  {
    \left|
      y_n-\hat{y}_n
    \right|
  }
  {
    \sqrt{1+\mathbf{z}'_n(\mathbfit{Z}'_{n-1}\mathbfit{Z}_{n-1})^{-1}\mathbf{z}_n}
    \hat{\sigma}_{n-1}
  }
\end{multline}
(cf.~(\ref{eq:pivotGLM});
the goodness of the definition follows from the formulas
given at the beginning of this section).
The expression on the right-hand side of (\ref{eq:nonconformityGLM})
can be replaced by other natural expressions,
such as $\lvert y_n-\hat y_n\rvert$.
See \cite{vovk/etal:2005book}, Section 8.5,
for further details.

According to our general convention,
the conformal predictor (\ref{eq:GammaGLM})
is called the \emph{Gauss predictor}
(although its discoverer was Fisher rather than Gauss).

We have already mentioned that the classical confidence predictor,
$\Gamma_n^{\epsilon}$ given by (\ref{eq:GammaGLM}),
does not work when there are many parameters;
in particular, it is required that $n\ge K+3$.
Theorem \ref{thm:2} shows that
there is hardly any way to use the knowledge that the first 10 explanatory variables
are the important ones
without abandoning the Gauss linear model:
no weakly valid confidence predictor in a very wide and natural class
can produce informative prediction intervals
unless $n\ge K+3$.
Indeed,
since the conditional distribution of the first $n$ observations
given $S_n$
is concentrated at one point for $n\le K+1$
and at two points for $n=K+2$ with probability one,
no conformal predictor
and, therefore, no weakly valid invariant confidence predictor
can give a bounded prediction region
$\Gamma_n^{\epsilon}$
for $\epsilon<0.5$ and $n\le K+2$.

\begin{remark}
  A common reaction to the importance of the condition $n\ge K+3$
  is that one can use only a subset of explanatory variables
  when $n<K+3$.
  A simple answer is that we are interested in confidence predictors
  that are valid under the Gauss linear model (\ref{eq:GL1}),
  not under some other model that is only ``approximately true'',
  in some ill-defined sense.
  However, the reader might wonder how the case of the Gauss linear model
  is different from the case of the IID model,
  where we are allowed to (and in Section \ref{sec:IID} do) take into account
  only a subset of explanatory variables before $n$ reaches $K+3$.
  The crucial difference is that
  if our data $(\mathbf{x}_n,y_n)$ conform to the IID model,
  discarding a subset of the explanatory variables
  will give us data $(\mathbf{x}'_n,y_n)$ that still conform to the IID model;
  on the other hand,
  if our data $(\mathbf{x}_n,y_n)$ conform to the Gauss linear model,
  discarding a subset of the explanatory variables
  will give us data $(\mathbf{x}'_n,y_n)$ that do not conform to the IID model,
  unless the coefficients $\beta_k$ in front of the discarded explanatory variables $x_{n,k}$
  (see (\ref{eq:GL1}))
  happen to be precisely zero.
  More generally,
  applying any transformation $\mathbf{x}'_n:=\phi(\mathbf{x}_n)$
  to the explanatory vectors does not lead us outside the IID model.
  The MVA model is similar to the IID model in this respect,
  except that the transformation $\phi$ has to be linear.
\end{remark}

\begin{figure}[bt]
  \centering
    \ifCOLOR
      \makebox{\includegraphics[width=\picturewidth]{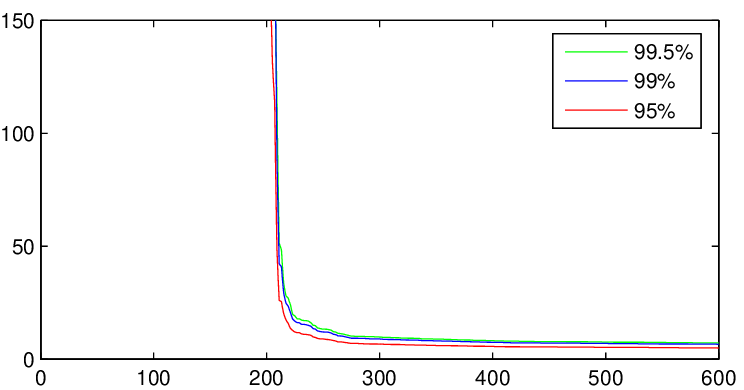}}
    \fi
    \ifnotCOLOR
      \makebox{\includegraphics[width=\picturewidth]{SM0_Med_bw.eps}}
    \fi
  \caption{\label{fig:SM-med}The median-accuracy plot for the classical prediction intervals.}
\end{figure}

Figure \ref{fig:SM-med} gives the median-accuracy plot
for the confidence predictor (\ref{eq:GammaGLM});
the predictor works very well soon after the number of observations reaches $K+3=103$.
Since the median is plotted,
the good quality of the prediction intervals shows only from $n=205$:
indeed, for $n<205$ at least half of the observed prediction intervals
are infinitely wide.

\section{The MVA model}
\label{sec:MVA}

Remember that the MVA model assumes, besides (\ref{eq:GL1}),
that $\mathbf{x}_n$ are generated independently
from the same unknown multivariate Gaussian distribution on $\bbbr^K$,
with the noise random variables $\xi_1,\xi_2,\ldots$
independent of $\mathbf{x}_1,\mathbf{x}_2,\ldots$\,.
The ATTS statistics in the MVA model are
\begin{equation*}
  S_n
  :=
  \left(
    \sum_{i=1}^n \mathbf{x}_i,
    \sum_{i=1}^n y_i,
    \sum_{i=1}^n \mathbf{x}_i\mathbf{x}'_i,
    \sum_{i=1}^n y_i\mathbf{x}_i,
    \sum_{i=1}^n y^2_i
  \right);
\end{equation*}
equivalently,
the ATTS statistics can be defined to be the empirical means and covariances of all variables,
i.e., the response and the explanatory variables.

Let $\mathbf{y}:=\mathbf{y}_n$,
$\mathbfit{Z}:=\mathbfit{Z}_n$,
$K^{\dagger}:=K^{\dagger}_n$
and $\mathbfit{U}$
be as in Section \ref{sec:IID}.
Suppose the value of the statistic $S_n$ is known.
The vector of residuals (\ref{eq:residuals1}) can now be written as
\begin{equation}\label{eq:residuals2}
  \mathbf{e}
  :=
  \mathbf{y}
  -
  \mathbfit{U}
  \left(
    \mathbfit{U}'\mathbfit{U}+a\mathbfit{I}
  \right)^{-1}
  \mathbfit{U}'\mathbf{y}
  =
  \mathbf{y}
  -
  \mathbfit{U}\mathbf{c},
\end{equation}
where
$\mathbf{c}:=(\mathbfit{U}'\mathbfit{U}+a\mathbfit{I})^{-1}\mathbfit{U}'\mathbf{y}$
is a known vector.
Since the joint distribution of $\mathbf{y}$ and the non-dummy columns of $\mathbfit{U}$
is invariant with respect to rotations around the vector $\boldsymbol{1}$,
the distribution of $\mathbf{e}$ will also be invariant
with respect to such rotations.
It might help the reader's intuition
to notice that knowing the value of $S_n$ is equivalent
to knowing the lengths of and the angles between
the following $K+2$ vectors:
the $K+1$ columns of $\mathbfit{Z}$ and $\mathbf{y}$.

In the rest of this section we will assume $n\ge3$
(with arbitrary conventions for $n=1,2$).
Let $e_1,\ldots,e_n$ be the components
of the vector (\ref{eq:residuals2}) of residuals
and $\overline{e}_{n-1}$ be the average of $e_1,\ldots,e_{n-1}$.
A standard statistical result
(stated in Section \ref{sec:on-lineII}; see (\ref{eq:pivotGM}))
allows us to conclude that
\begin{equation}\label{eq:pivotMVA}
  \sqrt{\frac{n-1}{n}}
  \frac
  {e_n-\overline{e}_{n-1}}
  {
    \sqrt
    {
      \frac{1}{n-2}
      \sum_{i=1}^{n-1}
      (e_i-\overline{e}_{n-1})^2
    }
  }
\end{equation}
has the $t$-distribution with $n-2$ degrees of freedom.

Let us see how to implement the conformal predictor corresponding
to the nonconformity measure
\begin{equation}\label{eq:nonconformityMVA}
  A_n
  \left(
    S_{n-1}(\mathbf{x}_1,y_1,\ldots,\mathbf{x}_{n-1},y_{n-1}),
    (\mathbf{x}_n,y_n)
  \right)
  :=
  \frac
  {e_n-\overline{e}_{n-1}}
  {
    \sqrt
    {
      \sum_{i=1}^{n-1}
      (e_i-\overline{e}_{n-1})^2
    }
  },
\end{equation}
which is proportional to (\ref{eq:pivotMVA});
the fact that the right-hand side of (\ref{eq:nonconformityMVA})
depends on the first $n-1$ observations only through the value of $S_{n-1}$
can be seen from the representation (\ref{eq:residuals2}),
where $\mathbf{c}$ is a known vector.
First we replace the true value $y_n$ by variable $y$ ranging over $\bbbr$.
Each residual $e_i$ becomes a linear
(according to (\ref{eq:residuals2}),
where $\mathbf{c}$ also depends on $y$)
function $e_i(y)$ of $y$,
and the prediction region can be written as
\begin{equation*}
  \Gamma_n^{\epsilon}
  :=
  \left\{
    y\in\bbbr
    \st
    \sqrt{\frac{n-1}{n}}
    \frac
    {
      \left|
        e_n(y)-\overline{e}_{n-1}(y)
      \right|
    }
    {
      \sqrt
      {
        \frac{1}{n-2}
        \sum_{i=1}^{n-1}
        (e_i(y)-\overline{e}_{n-1}(y))^2
      }
    }
    <
    t_{n-2}^{\epsilon/2}
  \right\}.
\end{equation*}
The inequality in this formula is quadratic in $y$,
so $\Gamma_n^{\epsilon}$ is easy to find.
We can see that the prediction region for $y_n$
is an interval (empirically, this is the typical case),
the union of two rays,
the empty set, or the whole real line.

Replacing $\Gamma_n^{\epsilon}$ by $\co\Gamma_n^{\epsilon}$
in the conformal predictor we have just defined
gives the \emph{MVA predictor}.
Our experiments with the artificial data set of Section \ref{sec:data}
are carried out as before
(cf.\ (\ref{eq:K})):
$\mathbfit{U}$ is defined as the first $11$ columns of $\mathbfit{Z}$ if $n<103$
and as the full $\mathbfit{Z}$ otherwise.

\begin{figure}[bt]
  \centering
    \ifCOLOR
      \makebox{\includegraphics[width=\picturewidth]{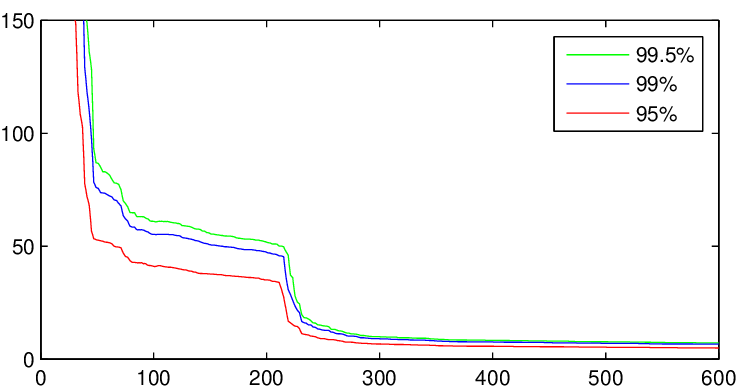}}
    \fi
    \ifnotCOLOR
      \makebox{\includegraphics[width=\picturewidth]{MVA0_Med_bw.eps}}
    \fi
  \caption{\label{fig:MVA-med}The median-accuracy plot for the MVA predictor.}
\end{figure}

The median-accuracy plot for the MVA predictor
and our artificial data set
is shown in Figure \ref{fig:MVA-med}.
Before the threshold $103$ the predictor quickly learns
$\alpha$ and the first $10$ parameters $\beta_k$,
and its performance more or less stabilizes
before quickly improving again when it starts learning the other parameters
from $n=103$ onwards;
the second improvement in the performance shows on the median-accuracy plot
from $n=205$.

The performance of the MVA predictor is better
than the performance of any other confidence predictor considered in this paper.
Of course,
this should not be taken to mean that the other predictors are worse.
Different predictors are based on different information about the data set.
None of the predictors ``knows'' that the components of $\mathbf{x}_n$
are realizations of independent standard Gaussian random variables;
even the MVA model,
the narrowest model considered in this paper,
allows arbitrary means of and arbitrary correlations between different explanatory variables
for the same observation.
The Gauss predictor does not know that the $\mathbf{x}_n$ are IID and Gaussian.
The IID predictor only knows that the observations $(\mathbf{x}_n,y_n)$ are IID,
and the IID--Gauss predictor,
introduced in the next section,
knows, in addition, that the $y_n$ are generated by (\ref{eq:GL1}).

The median-accuracy plot for each of the four predictors
is essentially determined by that for the MVA predictor
and the threshold for the corresponding model
as shown in Table \ref{tab:limitations}.
It is convenient to represent
each line on a median-accuracy plot
as the function that maps each value for the accuracy in the interval $[0,150]$
to the first step at which that accuracy is achieved
(so the graph of this function is obtained by rotating the page
by $90^{\circ}$ counterclockwise).
Each of the three functions in Figure \ref{fig:RRCM-med}
is, approximately, the maximum of $2\lceil 1/\epsilon\rceil$
and the corresponding function in Figure \ref{fig:MVA-med}.
Similarly,
each of the three functions in Figure \ref{fig:SM-med}
is, approximately, the maximum of $2(K+3)=206$
and the corresponding function in Figure \ref{fig:MVA-med}.
As usual, the factor of 2 appears because of the use of median in our accuracy plots.

\section{The IID--Gauss model}
\label{sec:IID-Gauss}

As defined in Section \ref{sec:introduction},
the IID--Gauss model is the combination
of the Gauss linear and IID models:
we assume both that the observations are IID
and that the responses are generated by (\ref{eq:GL1})
with $\xi_1,\xi_2,\ldots$ independent of $\mathbf{x}_1,\mathbf{x}_2,\ldots$\,.
Correspondingly, the ATTS statistics are
\begin{equation*}
  S_n
  :=
  \left(
    \lbag \mathbf{x}_1,\ldots,\mathbf{x}_n\rbag,
    \sum_{i=1}^n y_i,
    \sum_{i=1}^n y_i\mathbf{x}_i,
    \sum_{i=1}^n y^2_i
  \right).
\end{equation*}

\begin{figure}[bt]
  \centering
    \ifCOLOR
      \makebox{\includegraphics[width=\picturewidth]{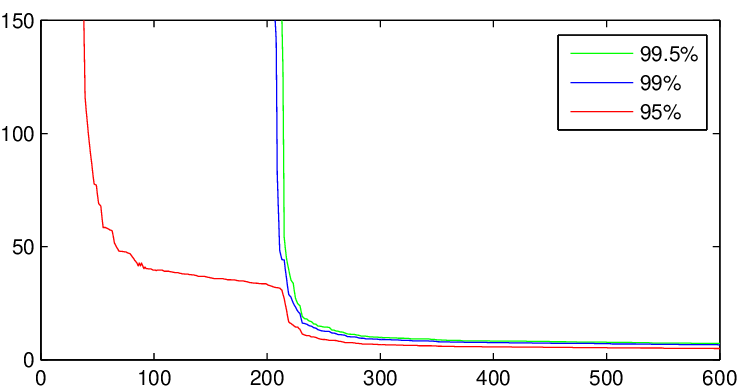}}
    \fi
    \ifnotCOLOR
      \makebox{\includegraphics[width=\picturewidth]{comb0_Med_bw.eps}}
    \fi
  \caption{\label{fig:comb-med}The median-accuracy plot
    for the IID--Gauss predictor.}
\end{figure}

Using the nonconformity measure (\ref{eq:nonconformity-de-Finetti})
and replacing the prediction regions
output by the corresponding conformal predictor
with their convex hulls,
we obtain the \emph{IID--Gauss predictor}.
Its performance on our usual data set is shown
in Figure \ref{fig:comb-med}.
We do not know whether the IID--Gauss predictor can be implemented efficiently,
and Figure \ref{fig:comb-med} was produced
using Monte-Carlo sampling from the conditional distributions given $S_n$.
However,
comparing Figure \ref{fig:comb-med}
to Figures \ref{fig:RRCM-med} (to the left of $n=205$)
and \ref{fig:SM-med} (to the right of $n=205$),
we can see that the following simple confidence predictor will work
almost as well as the IID--Gauss predictor
on our data set:
predict using the IID predictor if $n<103$
and predict using the Gauss predictor if $n\ge103$.
As in all other cases in this paper where the threshold $n=K+3=103$ appears,
the best switch-over point will be slightly greater than $K+3$,
but the question of when exactly to switch
is outside the scope of this paper.

\begin{remark}
  The IID predictor and the IID--Gauss predictor
  use the same nonconformity measure,
  (\ref{eq:nonconformity-de-Finetti}),
  but still produce very different median-accuracy plots
  at confidence level $99.5\%$.
  This happens because of the conditioning on the event $S_n\rnd=S_n\obs$
  in the definition (\ref{eq:p}).
  Since the ATTS statistics perform more radical data compression
  in the case of the IID--Gauss model,
  the achievable values of $\Prob(A_n\rnd\ge A_n\obs\given S_n\rnd=S_n\obs)$
  (corresponding to (\ref{eq:p}) with $\tau_n:=1$)
  are much smaller than the $1/n$ achievable under the IID model.
\end{remark}

As in the previous section,
there is a close connection between Figures \ref{fig:MVA-med} and \ref{fig:comb-med}:
each of the three functions in Figure \ref{fig:comb-med}
is, approximately, the maximum of $2\min(\lceil 1/\epsilon\rceil,K+3)$
and the corresponding function in Figure \ref{fig:MVA-med}.
The distributive law of $\max$ over $\min$
now implies that each of the three functions in Figure \ref{fig:comb-med}
is the minimum of the corresponding functions
in Figures \ref{fig:RRCM-med} and \ref{fig:SM-med}.

\section{On-line protocol, part II}
\label{sec:on-lineII}

Theorem~\ref{thm:1} sheds new light
not only on the main topic of this paper, predictive linear regression,
but also on some more classical corners of statistics.
In this section we will discuss, in particular,
Fisher's fiducial prediction and Wilks's non-parametric prediction intervals.
At the end of the section we discuss relaxations of the on-line protocol.

\subsection*{The Gaussian model}


Let us consider the model (\ref{eq:GL1}) with the $\mathbf{x}_n$ absent
(i.e., $K=0$);
in other words, $y_n$ is an IID sequence with $y_n\sim N(\alpha,\sigma^2)$
and unknown $\alpha$ and $\sigma^2>0$.
This model will be called the \emph{Gaussian model}.
Notice that the MVA model and the IID--Gauss model
also reduce to the Gaussian model when $K=0$.

The fact that
\begin{equation}\label{eq:pivotGM}
  T_n
  :=
  \sqrt{\frac{n-1}{n}}
  \frac{y_n-\overline{y}_{n-1}}{\hat{\sigma}_{n-1}},
\end{equation}
where
\begin{equation*}
  \overline{y}_l
  :=
  \frac{1}{l}
  \sum_{i=1}^l
  y_i
  \quad\text{and}\quad
  \hat{\sigma}^2_l
  :=
  \frac{1}{l-1}
  \sum_{i=1}^l
  (y_i-\overline{y}_l)^2,
\end{equation*}
has the $t$-distribution with $n-2$ degrees of freedom
\cite{fisher:1925applications}
allows us to conclude that
$y_n\in\Gamma^{\epsilon}_n$ with probability $1-\epsilon$,
where the prediction interval $\Gamma^{\epsilon}_n$ for $y_n$ is defined by
\begin{equation}\label{eq:GammaGM}
  \Gamma_n^{\epsilon}
  :=
  \left\{
    y\in\bbbr
    \st
    \left|
      y-\overline{y}_{n-1}
    \right|
    <
    t_{n-2}^{\epsilon/2}
    \sqrt{\frac{n}{n-1}}
    \hat{\sigma}_{n-1}
  \right\},
  \quad
  n=3,4,\ldots,
\end{equation}
and $\epsilon\in(0,1)$ is the chosen significance level.
This prediction interval is a special case of (\ref{eq:GammaGLM}).

Fisher discussed (\ref{eq:GammaGM}) and related confidence predictors
in his last book
(\cite{fisher:1973}, Sections V.3--4)
under the rubric of ``fiducial prediction''.
It appears that the idea of fiducial prediction
is less controversial
(and less often discussed)
than the related idea of fiducial inference for parameter values;
besides,
we will be interested in the least controversial aspects of fiducial prediction.
Fisher's comments about fiducial prediction in Sections V.3--4
are all applicable to the predictor (\ref{eq:GammaGM}),
although in Section V.3 he discusses prediction
of exponentially distributed rather than Gaussian random variables.

To some extent answering his critics
(``some teachers assert that statements of fiducial probability
cannot be tested by observations''),
he writes that ``fiducial statements about future observations''
(such as (\ref{eq:GammaGM}), although this passage
is about exponentially distributed responses)
``are verifiable by subsequent observations to any degree of precision required''.
The following is our reconstruction
(we believe the only possible reconstruction)
of \emph{Fisher's verification protocol},
as applied to the prediction intervals (\ref{eq:GammaGM}).
Fix a significance level $\epsilon\in(0,1)$ and $l\in\{2,3,\ldots\}$
(the \emph{sample size};
we might consider samples of different sizes,
but we will stick to the simplest case).
For $m=1,2,\ldots$,
generate the \emph{$m$th sample}
\begin{equation*}
  y_{(m-1)(l+1)+1},
  y_{(m-1)(l+1)+2},
  \ldots,
  y_{m(l+1)-1}
\end{equation*}
and the \emph{$m$th test observation} $y_{m(l+1)}$.
Register an error if the $m$th prediction interval
computed from the $m$th sample according to (\ref{eq:GammaGM})
fails to contain the $m$th test observation:
\begin{equation*}
  \err^{\dagger}_m
  :=
  \begin{cases}
    0 & \text{if $\lvert y_{m(l+1)}-\overline{y}\rvert
                 < t^{\epsilon/2}_{l-1}\sqrt{\frac{l+1}{l}}
                 \sqrt{\frac{1}{l-1}\sum_{i=(m-1)(l+1)+1}^{m(l+1)-1}(y_i-\overline{y})^2}$}\\
    1 & \text{otherwise},
  \end{cases}
\end{equation*}
where
\begin{equation*}
  \overline{y}
  :=
  \frac{1}{l}
  \sum_{i=(m-1)(l+1)+1}^{m(l+1)-1}
  y_{i}.
\end{equation*}
As in the on-line protocol,
the errors $\err^{\dagger}_m$, $m=1,2,\ldots$, are independent.
The frequency of error gets arbitrarily close to $\epsilon$
with an arbitrarily high probability
as the number of observations increases.

Fisher's verification protocol has a serious drawback:
as Fisher puts it,
\begin{quotation}
  In carrying out such a verification [\dots],
  it is to be supposed that the investigator
  is not deflected from his purpose
  by the fact that new data are becoming available
  from which predictions,
  better than the one he is testing,
  could at any time be made.
  For verification,
  the original prediction must be held firmly in view.
  This, of course,
  is a somewhat unnatural attitude for a worker
  whose main preoccupation is to improve his ideas.
\end{quotation}
Indeed,
when making his prediction for the $m$th test observation,
the ``investigator'' is asked to ignore the first $m-1$ samples.
The protocol seems to be an artificial device
rather than a description of what
``a worker whose main preoccupation is to improve his ideas''
might do in reality.
Let us see, however, what happens if all the previous observations
\emph{are} used when making the $m$th prediction;
in this case,
the sequence of errors becomes
\begin{equation*}
  \!\!\!\!\!
  \err^{\ddagger}_m
  :=
  \begin{cases}
    0 & \text{if $\lvert y_{m(l+1)}-\overline{y}\rvert
                 < t^{\epsilon/2}_{m(l+1)-2}\sqrt{\frac{m(l+1)}{m(l+1)-1}}
                 \sqrt{\frac{1}{m(l+1)-2}\sum_{i=1}^{m(l+1)-1}(y_{i}-\overline{y})^2}$}\\
    1 & \text{otherwise},
  \end{cases}
\end{equation*}
where
\begin{equation*}
  \overline{y}
  :=
  \frac{1}{m(l+1)-1}
  \sum_{i=1}^{m(l+1)-1} y_{i}.
\end{equation*}
As $\err^{\ddagger}_m$, $m=1,2,\ldots$,
is a subsequence of the sequence of errors $\err^{\epsilon}_n$, $n=1,2,\ldots$,
in the on-line protocol,
the errors are still independent.
Theorem~\ref{thm:1} cures the drawback.

Fisher's theory of fiducial prediction is based on the fact
that a value such as (\ref{eq:pivotGM}) has a known distribution for each $n$;
therefore, it can be used as a ``pivot'' to project this known distribution
onto the future observation $y_n$.
This idea may be difficult to formalize,
but Fisher's observation that (\ref{eq:pivotGM}) has a known distribution
can be strengthened:
Theorem~\ref{thm:1}
(applied to the nonconformity measure (\ref{eq:pivotGM}))
implies that the random variables $T_n$, $n=3,4,\ldots$,
have the $t$-distribution with $n-2$ degrees of freedom
and are independent in the on-line protocol.
Therefore, not only the individual $T_n$ have known distributions,
but also the whole sequence $(T_1,T_2,\ldots)$ has a known distribution
(the product of $t$-distributions).

\subsection*{The univariate IID model}

The IID model is different from all the other models in this paper
(see Figure \ref{fig:models})
in that it gives a univariate model different from the Gaussian model
in the case where the explanatory variables are absent.
The construction of prediction and tolerance intervals
in the univariate IID model,
which says that $y_1,y_2,\ldots$ form an IID sequence,
was undertaken by many authors
following the pioneering paper by Wilks \cite{wilks:1941}.
Wilks's work was later extended to the multivariate case:
see, e.g., Fraser \cite{fraser:1957};
this extension, however,
is not directly related to our IID predictors.
For simplicity,
let us assume in this subsection,
as is customary in literature,
that the distribution of one observation is continuous.
Correspondingly,
we will assume that the realized values of $y_n$, $n=1,2,\ldots$,
are all different.

For each $n=1,2,\ldots$,
define $T_n\in\{1,2,\ldots,n\}$ as the smallest $i$ such that
$y_n<y_{(n-1,i)}$,
where $y_{(n-1,1)},\ldots,y_{(n-1,n-1)}$ is the sequence of the first $n-1$ observations
$y_1,\ldots,y_{n-1}$ sorted in the ascending order;
if $y_n>y_{(n-1,n-1)}$,
set $T_n:=n$.
Each $T_n$ is a ``pivot'',
being distributed uniformly on the set $\{1,\ldots,n\}$.
Wilks suggested the following prediction intervals based on this fact:
fix a number $r\in\{1,2,\ldots\}$ and define
$\Gamma_n^{2r/n}$, $n=2r+1,2r+2,\ldots$,
to be the interval $(y_{(n-1,r)},y_{(n-1,n-r)})$;
the probability of error, $y_n\notin\Gamma_n^{2r/n}$,
is then $2r/n$.
Now Theorem~\ref{thm:1} implies that the whole random sequence $(T_1,T_2,\ldots)$
has a known distribution:
namely, it is distributed according to the product $U_1\times U_2\times\cdots$
of the uniform distributions $U_n$ on $\{1,\ldots,n\}$.
In particular, Wilks's prediction intervals
$\Gamma_n^{2r/n}$, $n=2r+1,2r+2,\ldots$,
lead to independent errors.

\subsection*{Relaxations of the on-line protocol}


This paper concentrates on the on-line prediction protocol.
Smoothed conformal predictors lead to independent errors
in the on-line protocol,
and Theorem~\ref{thm:2} suggests that
conformal predictors are the most natural
weakly valid confidence predictors.
This is why we included the requirement of independence
in the definition of strong validity,
despite the fact that the error frequency can be shown to approach
the error probability $\epsilon$
with probability approaching one
even when the requirement of independence is relaxed in certain ways.


The situation changes when we move outside the on-line protocol.
The on-line protocol is natural,
but in one respect it is overly restrictive:
the true response $y_n$ becomes known
before the prediction for the next response $y_{n+1}$ is made.
It can be shown that the error frequency will still converge to $\epsilon$
if the true response is only given for a small fraction of observations,
and even for those observations it can be given with a delay
(\cite{vovk/etal:2005book}, Section 4.3;
see also \cite{vanderlooy/etal:2007} for an empirical study).
The independence of errors, however, will be lost
(we can still have ``approximate independence'',
but this is a much more elusive notion
than ordinary independence).

\section{Conclusion}
\label{sec:conclusion}

In this paper we considered the problem of prediction
in three main regression models.
One of these models,
the Gauss linear model,
is the standard textbook one.
The MVA model seems to have been somewhat neglected,
partly because of philosophical reasons:
according to the conditionality principle
(\cite{cox/hinkley:1974}, Section 2.3(iii))
one should condition on the observed values of the explanatory variables
to make the prediction (or estimate, etc.)\ more relevant
to the data at hand.
In most of this paper we took a pragmatic approach,
studying which models permit one to produce
informative prediction intervals in different circumstances
without being restricted \emph{a priori} by general principles.
We did use the sufficiency principle
in our interpretation of Theorem~\ref{thm:2},
but we admit that this makes the theorem less convincing.
Surprisingly,
the IID model appears to have been 
neglected in the field of regression,
even in non-parametric statistics,
where the value of this model is in principle well understood.

\subsection*{Acknowledgments}

  This paper has greatly benefited from Glenn Shafer's advice,
  from a discussion with Steffen Lauritzen,
  and from correspondence with Philip Dawid.
  Comments by the anonymous referees of the journal version
  and by Prof.\ Susan Murphy
  helped us improve the presentation.
  This work was partially supported by EPSRC (grant EP/F002998/1),
  MRC (grant S505/65),
  and the Royal Society.

  \ifarXiv
    \appendix
    \section*{Appendix~A: Proofs of the theorems}
  \fi

  In this appendix we will prove
  the two main results stated in this paper,
  Theorems~\ref{thm:1} and~\ref{thm:2}.
  A version of Theorem \ref{thm:1} was proved in Section 8.7 of \cite{vovk/etal:2005book},
  but we reproduce the principal points of the proof
  to make our exposition self-contained.
  A special case of Theorem~\ref{thm:2}
  (namely, for the IID model)
  was proved in Section 2.6 of \cite{vovk/etal:2005book}.

  \subsection*{Proof of Theorem~\ref{thm:1}}

  In this proof,
  $\zeta_1,\zeta_2,\ldots$
  will be random observations generated by $P\in\PPP$,
  $(\zeta_1,\zeta_2,\ldots)\sim P$,
  and $\tau_1,\tau_2,\ldots$ will be random numbers,
  $(\tau_1,\tau_2,\ldots)\sim U^{\infty}$.
  For each $n=0,1,\ldots$ let $\GGG_{n}$ be the $\sigma$-algebra
  generated by the random elements
  \begin{equation*}
    S_n(\zeta_1,\ldots,\zeta_n),
    \zeta_{n+1},\tau_{n+1},
    \zeta_{n+2},\tau_{n+2},\ldots\,.
  \end{equation*}
  So $\GGG_0$ is the most informative $\sigma$-algebra
  and $\GGG_0\supseteq\GGG_1\supseteq\GGG_2\supseteq\cdots$.
  It will be convenient to write $\Prob_{\GGG}(E)$ and $\Expect_{\GGG}(\xi)$
  for the conditional probability $\Prob(E\given\GGG)$
  and expectation $\Expect(\xi\given\GGG)$,
  respectively,
  given a $\sigma$-algebra $\GGG$.

  \begin{lemma}\label{lem:basic}
    For any step $n=1,2,\ldots$ and any $\epsilon\in(0,1)$,
    \begin{equation*}
      \Prob_{\GGG_n}
      \left(
        p_n\le\epsilon
      \right)
      =
      \epsilon.
    \end{equation*}
  \end{lemma}
  \begin{proof}
    For a given value of the summary $S_n(\zeta_1,\ldots,\zeta_n)$
    of the first $n$ observations,
    consider the conditional distribution function $F$ of the random variable
    $\eta:=A_n(S_{n-1}(\zeta_1,\ldots,\zeta_{n-1}),\zeta_n)$
    (because of the total sufficiency,
    it does not matter whether we further condition on
    $\zeta_{n+1},\tau_{n+1},\zeta_{n+2},\tau_{n+2},\ldots$).
    Define $F(x-)$ to be $\sup_{t<x}F(t)$.
    Our task is to show that the conditional probability of the event
    \begin{equation}\label{eq:goal}
      1-F(\eta)+\tau_n(F(\eta)-F(\eta-))\le\epsilon
    \end{equation}
    is $\epsilon$
    (since the left-hand side of (\ref{eq:goal}) coincides
    with the right-hand side of the definition (\ref{eq:p})).
    The latter fact is usually stated in statistics textbooks
    for continuous $F$
    (see, e.g., \cite{cox/hinkley:1974}, p.~66),
    but it is also easy to check in general.
  \end{proof}
  \begin{lemma}\label{lem:measurable}
    For any step $n=1,2,\ldots$,
    $p_n$ is $\GGG_{n-1}$-measurable.
  \end{lemma}
  \begin{proof}
    This follows from the definition:
    $p_n$ is defined in terms of $\zeta_n$, $\tau_n$
    and the summary of the first $n-1$ observations.
  \end{proof}

  Now we can easily prove the theorem.
  First we demonstrate that,
  for any $n=1,2,\ldots$ and any $\epsilon_1,\ldots,\epsilon_n\in(0,1)$,
  \begin{equation}\label{eq:conditional}
    \Prob_{\GGG_n}\!
    \left(
      p_n\le\epsilon_n,\ldots,p_1\le\epsilon_1
    \right)
    =
    \epsilon_n\cdots\epsilon_1
    \quad
    \text{a.s.}
  \end{equation}
  The proof is by induction on $n$.
  For $n=1$,
  (\ref{eq:conditional}) is a special case of Lemma~\ref{lem:basic}.
  For $n>1$ we obtain,
  from Lemmas~\ref{lem:basic} and~\ref{lem:measurable},
  standard properties of conditional expectations,
  and the inductive assumption:
  \begin{multline*}
    \Prob_{\GGG_n}\!
    \left(
      p_n\le\epsilon_n,\ldots,p_1\le\epsilon_1
    \right)
    =
    \Expect_{\GGG_n}\!
    \left(
      \Expect_{\GGG_{n-1}}
      \left(
        \III_{p_n\le\epsilon_n}
        \III_{p_{n-1}\le\epsilon_{n-1},\ldots,p_1\le\epsilon_1}
      \right)
    \right)
  \\
    {}=
    \Expect_{\GGG_n}\!
    \left(
      \III_{p_n\le\epsilon_n}
      \Expect_{\GGG_{n-1}}
      \left(
        \III_{p_{n-1}\le\epsilon_{n-1},\ldots,p_1\le\epsilon_1}
      \right)
    \right)
    =
    \Expect_{\GGG_n}\!
    \left(
      \III_{p_n\le\epsilon_n}
      \epsilon_{n-1}\cdots\epsilon_1
    \right)
  \\
    =
    \epsilon_n\epsilon_{n-1}\cdots\epsilon_1
    \quad
    \text{a.s.}
  \end{multline*}
  The ``tower property'' of conditional expectations
  immediately implies
  \begin{equation*}
    \Prob
    \left(
      p_n\le\epsilon_n,\ldots,p_1\le\epsilon_1
    \right)
    =
    \epsilon_n\cdots\epsilon_1.
  \end{equation*}
  Therefore, the distribution of the first $n$ $p$-values $p_1,\ldots,p_n$
  is $U^n$, for all $n=1,2,\ldots$\,.
  This implies that the distribution of the infinite sequence $p_1,p_2,\ldots$
  is $U^{\infty}$.

  \subsection*{Proof of Theorem~\ref{thm:2}}


  In this proof, $\Zeta:=\bbbr^K\times\bbbr$
  and $\zeta_i$ stands for $(\mathbf{x}_i,y_i)$.
  Let 
  $n\in\NNN$.

  For each summary $s\in\Sigma_n$
  let $f(s)$ be the conditional probability
  given $S_n(\zeta_1,\ldots,\zeta_n)=s$
  that $\Gamma$ makes an error at a significance level $\epsilon$
  when predicting $y_n$ from $\zeta_1,\ldots,\zeta_{n-1}$ and $\mathbf{x}_n$,
  the observations $\zeta_1,\zeta_2,\ldots$ being generated from $P\in\PPP$.
  We know that the expected value of $f(S_n(\zeta_1,\ldots,\zeta_n))$ is $\epsilon$
  under any $P\in\PPP$, and this,
  by the bounded completeness of $S_n$,
  implies that $f(s)=\epsilon$
  for almost all (under $PS_n^{-1}$ for any $P\in\PPP$) summaries $s$.
  Define $E(s,\epsilon)$ to be the set of all pairs
  $(s',\zeta)=(s',(\mathbf{x},y))\in\Sigma_{n-1}\times\Zeta$
  such that $F_n(s',\zeta)=s$
  (where $F_n$ is the function from the definition
  of the algebraic transitivity of the $S_n$)
  and $\Gamma$ makes an error at the significance level $\epsilon$
  when predicting $y$ and fed with $\zeta_1,\ldots,\zeta_{n-1}$
  satisfying $S_{n-1}(\zeta_1,\ldots,\zeta_{n-1})=s'$ and with $\mathbf{x}$
  (since $\Gamma$ is invariant, whether an error is made depends only
  on $s'$, not on the particular $\zeta_1,\ldots,\zeta_{n-1}$).
  It is clear that
  \begin{equation*}
    \epsilon_1\le\epsilon_2
    \Longrightarrow
    E(s,\epsilon_1)
    \subseteq
    E(s,\epsilon_2)
  \end{equation*}
  and
  \begin{equation*}
    \Prob
    \left(
      (S_{n-1}(\zeta_1,\ldots,\zeta_{n-1}),\zeta_n)\in E(s,\epsilon)
      \given
      S_n(\zeta_1,\ldots,\zeta_n) = s
    \right)
    =
    \epsilon
    \quad
    \text{a.s.},
  \end{equation*}
  where $(\zeta_1,\zeta_2,\ldots)\sim P\in\PPP$.

  In this proof we say ``conformity measure''
  to mean a nonconformity measure
  which is used for computing $p$-values in the opposite way
  to (\ref{eq:p}):
  the ``$>$'' in (\ref{eq:p}) is replaced by ``$<$''.
  Let us check that the conformal predictor $\Gamma^{\dagger}$
  determined by the conformity measure
  \begin{equation*}
    A_n(s',\zeta)
    :=
    \inf
    \left\{
      \epsilon
      \st
      (s',\zeta)\in E(F_n(s',\zeta),\epsilon)
    \right\}
  \end{equation*}
  is at least as accurate as $\Gamma$.
  By the monotone convergence theorem for conditional expectations,
  \begin{multline*}
    \Prob
    \left(
      A_n(S_{n-1}(\zeta_1,\ldots,\zeta_{n-1}),\zeta_n) \le \epsilon
      \given
      S_n(\zeta_1,\ldots,\zeta_n) = s
    \right)
    \\
    {}=
    \lim_{\delta\downarrow\epsilon}
    \Prob
    \left(
      A_n(S_{n-1}(\zeta_1,\ldots,\zeta_{n-1}),\zeta_n) < \delta
      \given
      S_n(\zeta_1,\ldots,\zeta_n) = s
    \right)
    \\
    {}\le
    \lim_{\delta\downarrow\epsilon}
    \Prob
    \left(
      (S_{n-1}(\zeta_1,\ldots,\zeta_{n-1}),\zeta_n)\in E(s,\delta)
      \given
      S_n(\zeta_1,\ldots,\zeta_n) = s
    \right)\\
    =
    \lim_{\delta\downarrow\epsilon}\delta
    =
    \epsilon
    \quad
    \text{a.s.},
  \end{multline*}
  where $(\zeta_1,\zeta_2,\ldots)\sim P\in\PPP$
  and $\delta$ is constrained to be a rational number.
  Therefore,
  at each significance level $\epsilon$
  and for $P$-almost all $(\zeta_1,\ldots,\zeta_n)\in\Zeta^n$,
  \begin{multline*}
    y_n \in (\Gamma^{\dagger})^{\epsilon}(\zeta_1,\ldots,\zeta_{n-1},\mathbf{x}_n)
    \Longleftrightarrow
    \Prob
    \left(
      A_n\rnd
      \le
      A_n\obs
      \given
      S_n\rnd = S_n\obs
    \right)
    >
    \epsilon
    \\
    {}\Longrightarrow
    A_n\obs
    >
    \epsilon
    \Longrightarrow
    (S_{n-1}(\zeta_1,\ldots,\zeta_{n-1}),\zeta_n)
    \notin
    E(S_n(\zeta_1,\ldots,\zeta_n),\epsilon)\\
    \Longleftrightarrow
    y_n \in \Gamma^{\epsilon}(\zeta_1,\ldots,\zeta_{n-1},\mathbf{x}_n),
  \end{multline*}
  in the notation of (\ref{eq:p}).

  \ifarXiv
    \section*{Appendix~B: Explicit algorithms}
  \fi

  For all three main models considered in this paper
  there are relatively efficient algorithms for computing prediction intervals.
  In this appendix we will describe two of them,
  one new (the MVA predictor)
  and the other not known in the statistical community
  (the IID predictor).

  \subsection*{The MVA predictor}

  The idea of this algorithm was described in Section \ref{sec:MVA},
  and here we will give explicit formulas
  which the reader might find useful
  in writing computer programs implementing this idea.
  Let $\mathbfit{P}$ be the matrix
  \begin{equation*}
    \mathbfit{P}
    =
    \mathbfit{I}
    -
    \mathbfit{U}
    \left(
      \mathbfit{U}'\mathbfit{U}+a\mathbfit{I}
    \right)^{-1}
    \mathbfit{U}'
  \end{equation*}
  projecting the responses $\mathbf{y}$ onto the residuals $\mathbf{e}$
  (see (\ref{eq:residuals2})).
  Introducing the notation
  \begin{equation*}
    \mathbf{y}
    =
    \begin{pmatrix}
      y_1\\
      \vdots\\
      y_{n-1}\\
      y
    \end{pmatrix}
    =
    \begin{pmatrix}
      y_1\\
      \vdots\\
      y_{n-1}\\
      0
    \end{pmatrix}
    +
    \begin{pmatrix}
      0\\
      \vdots\\
      0\\
      y
    \end{pmatrix}
    =
    \mathbf{y}_0 + y\mathbf{u}
  \end{equation*}
  (so that $\mathbf{y}_0$ is the vector of the first $n-1$ responses followed by a $0$
  and $\mathbf{u}$ is the vector whose only non-zero component is a $1$
  at the $n$th position),
  we can define the vector of residuals as
  \begin{equation}\label{eq:e}
    \mathbf{e}
    =
    \mathbfit{P}(\mathbf{y}_0+y\mathbf{u})
    =
    \mathbfit{P}\mathbf{y}_0 + y\mathbfit{P}\mathbf{u}.
  \end{equation}
  Subtract from each component of $\mathbfit{P}\mathbf{y}_0$
  the mean of the first $n-1$ components of $\mathbfit{P}\mathbf{y}_0$
  and subtract from each component of $\mathbfit{P}\mathbf{u}$
  the mean of the first $n-1$ components of $\mathbfit{P}\mathbf{u}$;
  let the resulting vectors be $\mathbf{a}$ and $\mathbf{b}$,
  respectively.
  The prediction region consists of the $y$ satisfying
  \begin{equation}\label{eq:quadratic}
    \sqrt{\frac{n-1}{n}}
    \frac
    {
      \left|
        a_n + y b_n
      \right|
    }
    {
      \sqrt
      {
        \frac{1}{n-2}
        \sum_{i=1}^{n-1}
        (a_i+yb_i)^2
      }
    }
    <
    t,
  \end{equation}
  where $a_i$ and $b_i$ are the components of $\mathbf{a}$ and $\mathbf{b}$,
  respectively,
  and $t$ stands for $t_{n-2}^{\epsilon/2}$.
  We can rewrite (\ref{eq:quadratic}) as
  \begin{equation*}
    (a_n + b_n y)^2
    <
    t^2
    \frac{n}{(n-1)(n-2)}
    \sum_{i=1}^{n-1}
    (a_i + b_i y)^2,
  \end{equation*}
  or, explicitly as a quadratic inequality,
  \begin{equation*}
    A y^2 + 2 B y + C < 0,
  \end{equation*}
  where
  \begin{align*}
    A
    &=
    (n-1)(n-2) b_n^2
    -
    t^2 n \sum_{i=1}^{n-1} b_i^2,\\
    B
    &=
    (n-1)(n-2) a_n b_n
    -
    t^2 n \sum_{i=1}^{n-1} a_i b_i,\\
    C
    &=
    (n-1)(n-2) a_n^2
    -
    t^2 n \sum_{i=1}^{n-1} a_i^2.
  \end{align*}
  The solutions to the corresponding quadratic equation are
  \begin{equation}\label{eq:solutions}
    \frac
    {
      -B
      \pm
      \sqrt{D}
    }
    {
      A
    },
  \end{equation}
  where the discriminant $D$ is
  \begin{equation*}
    D
    =
    B^2-AC.
  \end{equation*}
  Therefore,
  the prediction interval
  (defined, as usual,
  to be the convex hull of the prediction region)
  is:
  \begin{itemize}
  \item
    unbounded if
    \begin{equation*}
      A < 0
    \end{equation*}
    or
    \begin{equation*}
      A = 0
      \text{ and }
      B \ne 0
    \end{equation*}
    or
    \begin{equation*}
      A = B = 0
      \text{ and }
      C < 0;
    \end{equation*}
  \item
    empty if
    \begin{equation*}
      A > 0
      \text{ and }
      D\le0
    \end{equation*}
    or
    \begin{equation*}
      A = B = 0
      \text{ and }
      C \ge 0;
    \end{equation*}
  \item
    otherwise (i.e., if $A>0$ and $D>0$),
    the open interval with the end-points (\ref{eq:solutions}).
  \end{itemize}

  \subsection*{The IID predictor}

  In this subsection we describe explicitly the IID predictor.
  We will follow mainly \cite{vovk/etal:2005book}, Section 2.3,
  although our current algorithm is simpler,
  as we are only interested in finding the prediction intervals
  (the convex hulls of the prediction regions)
  and not the prediction regions themselves.

  Rewrite the vector of residuals (\ref{eq:e}) in the form $\mathbf{a}+y\mathbf{b}$,
  where $\mathbf{a}:=\mathbfit{P}\mathbf{y}_0$ and $\mathbf{b}:=\mathbfit{P}\mathbf{u}$.
  For each $i=1,\ldots,n$, let
  \begin{equation*} 
    S_i
    :=
    \{
      y
      \st
      \lvert e_i(y)\rvert \ge \lvert e_n(y)\rvert
    \}
    =
    \{
      y
      \st
      \lvert a_i+b_i y\rvert \ge \lvert a_n + b_n y\rvert
    \},
  \end{equation*}
  where $a_i$ and $b_i$ are the components of $\mathbf{a}$ and $\mathbf{b}$.
  Each set $S_i$ (always closed) will be
  the real line,
  the union of two rays, a ray,
  an interval (by an interval we will always mean a bounded interval),
  a point (which is a special case of an interval), or empty.
  Indeed, as we are interested in $\lvert a_i + b_i y\rvert$ we can assume $b_i \ge 0$
  for $i=1,\ldots,n$
  (if necessary, multiply both $a_i$ and $b_i$ by $-1$).
  If $b_i \ne b_n$ then $\lvert e_i(y)\rvert$ and $\lvert e_n(y)\rvert$ are equal at two points
  (which may coincide):
  \begin{equation}\label{eq:two}
    -\frac{a_i - a_n}{b_i - b_n}
    \quad \text{and} \quad
    -\frac{a_i + a_n}{b_i + b_n};
  \end{equation}
  in this case, $S_i$ is an interval (possibly a point), the union of two rays,
  or the real line.
  If $b_i = b_n \ne 0$ then $\lvert e_i(y)\rvert=\lvert e_n(y)\rvert$ at the only point
  \begin{equation}\label{eq:one}
    -\frac{a_i + a_n}{2 b_i}
  \end{equation}
  (and $S_i$ is a ray)
  unless $a_i = a_n$, in which case $S_i = \bbbr$.
  If $b_i=b_n=0$, $S_i$ is either $\emptyset$ or $\bbbr$.

  To calculate the $p$-value $p(y)$ for any potential response $y$ of $x_n$,
  we count how many $S_i$ include $y$ and divide by $n$:
  \begin{equation*}
    p(y)
    =
    \frac{\left|\{i=1,\ldots,n\st y\in S_i\}\right|}{n}.
  \end{equation*}
  This formula immediately implies that each $\Gamma_n^{\epsilon}$,
  as defined by (\ref{eq:Gamma}),
  is a closed set:
  setting $k:=\lfloor\epsilon n\rfloor+1$,
  \begin{multline*}
    \Gamma_n^{\epsilon}
    =
    \left\{
      y
      \st
      p(y) > \epsilon
    \right\}
    =
    \left\{
      y
      \st
      \left|\{i=1,\ldots,n\st y\in S_i\}\right| \ge k
    \right\}\\
    =
    \bigcup_{\{i_1,\ldots,i_k\}\subseteq\{1,\ldots,n\}}
    S_{i_1} \cap \cdots \cap S_{i_k};
  \end{multline*}
  the last expression is a finite union of closed sets.

  The following algorithm implicitly keeps track of the number $M(y)$ of $i=1,\ldots,n-1$
  such that $y \in S_i$ for each $y\in\bbbr$.
  For each critical point $y_j$, $j=1,\ldots,m$,
  let $\NM(j)$ be the change in $M(y)$, as $y$ increases,
  at $y_j$ due to $y_j$ being a critical point.
  The set of critical points will be represented as a sequence $P$
  (initially empty);
  the critical points will be listed in $P$ and in $\NM$ in the same order.
  The algorithm will also compute the number $L$ of $S_i$
  that include points to the left of the left-most critical point
  and the number $R$ of $S_i$
  that include points to the right of the right-most critical point.
  It is clear that $L$ and $\NM$ (or $R$ and $\NM$) determine $M$.
  The word ``add'' in the description of the algorithm
  means ``attach at the end'' (unless we say explicitly ``add in front'')
  and has nothing to do with arithmetic addition.

  The algorithm is given a significance level $\epsilon$
  and outputs the corresponding prediction interval
  $\bar\Gamma_n^{\epsilon}:=\co\Gamma_n^{\epsilon}$.
  We will assume that the set of critical points
  ((\ref{eq:two}) and (\ref{eq:one}))
  is not empty
  and will supply the matrix $\mathbfit{U}$ with a lower index
  to explicitly indicate the dependence on $n$.

  \bigskip

  \noindent
  \textsc{IID predictor}

  \medskip

  \parshape=39
  \IndentI   \WidthI
  \IndentI   \WidthI
  \IndentI   \WidthI
  \IndentI   \WidthI	
  \IndentII  \WidthII
  \IndentI   \WidthI	
  \IndentI   \WidthI
  \IndentI   \WidthI	
  \IndentII  \WidthII	
  \IndentIII \WidthIII
  \IndentIII \WidthIII	
  \IndentIV  \WidthIV
  \IndentIII \WidthIII	
  \IndentIII \WidthIII	
  \IndentIV  \WidthIV
  \IndentIII \WidthIII	
  \IndentIII \WidthIII	
  \IndentIV  \WidthIV
  \IndentIII \WidthIII	
  \IndentII  \WidthII	
  \IndentIII \WidthIII
  \IndentII  \WidthII	
  \IndentIII \WidthIII
  \IndentIII \WidthIII	
  \IndentIV  \WidthIV
  \IndentIII \WidthIII	
  \IndentIV  \WidthIV
  \IndentIII \WidthIII	
  \IndentII  \WidthII	
  \IndentIII  \WidthIII
  \IndentII   \WidthII	
  \IndentI  \WidthI	
  \IndentI  \WidthI
  \IndentI  \WidthI
  \IndentI  \WidthI
  \IndentI  \WidthI
  \IndentI  \WidthI
  \IndentI  \WidthI
  \IndentI  \WidthI
  \noindent
    $\mathbfit{P} := \mathbfit{I} -
      \mathbfit{U}_n(\mathbfit{U}'_n\mathbfit{U}_n+a\mathbfit{I})^{-1}\mathbfit{U}'_n$;\\
    $\mathbf{a}=(a_1,\ldots,a_n)' := \mathbfit{P}(y_1,\ldots,y_{n-1},0)'$;\\
    $\mathbf{b}=(b_1,\ldots,b_n)' := \mathbfit{P}(0,\ldots,0,1)'$;\\
    FOR $i=1,\ldots,n$\\
      IF ($b_i<0$) $a_i:=-a_i$; $b_i:=-b_i$ END IF\\
    END FOR\\
    $P:=\NM:=()$; $L:=R:=0$;\\
    FOR $i=1,\ldots,n-1$\\
      IF ($b_i\ne b_n$)\\
         add the two points~(\ref{eq:two}) to $P$;\\
         IF (the corresponding $S_i$ is an interval)\\
           add $1$ and $-1$ to $\NM$\\
         END IF\\
         IF (the corresponding $S_i$ is a union of disjoint rays)\\
           add $-1$ and $1$ to $\NM$ and increase $L$ and $R$ by $1$\\
         END IF\\
         IF (the corresponding $S_i$ is the real line)\\
           reverse adding~(\ref{eq:two}) to $P$ and increase $L$ and $R$ by $1$\\
         END IF\\
      ELSEIF ($a_i=a_n$)\\
        increase $L$ and $R$ by $1$\\
      ELSEIF ($b_n\ne0$)\\
        add~(\ref{eq:one}) to $P$;\\
        IF (the corresponding $S_i$, which is a ray, is going to the right)\\
          add $1$ to $\NM$ and increase $R$ by 1\\
        ELSE\\
          add $-1$ to $\NM$ and increase $L$ by 1\\
        END IF\\
      ELSEIF ($\lvert a_i\rvert\ge\lvert a_n\rvert$)\\
        increase $L$ and $R$ by $1$\\
      END IF\\
    END FOR\\
    add $-\infty$ in front of $P$ and $L+1$ in front of $\NM$\\
    add $\infty$ at the end of $P$ and $-R-1$ at the end of $\NM$\\
    sort $P$ in ascending order breaking ties using $-\NM$;\\
    apply the same permutation to $\NM$;\\
    find the smallest $i_1$
      such that the sum of the first $i_1$ elements of $\NM$ exceeds $\epsilon$;\\
    find the largest $i_2$
      such that the sum of the first $i_2$ elements of $\NM$ exceeds $\epsilon$;\\
    $\bar\Gamma^{\epsilon}_n := [P(i_1),P(i_2)]$.

  \bigskip

  \noindent
  This algorithm is run by Predictor at each step $n=1,2,\ldots$
  of the on-line prediction protocol.
  It is important that both points~(\ref{eq:two}) should be added to $P$,
  even when they coincide.
  Finding $i_2$ can be done by search from the right,
  since the sum of all elements of $\NM$ is $0$.

  Let us suppose that the number $K$ of explanatory variables is bounded above by a constant.
  The computation time of our algorithm is $O(n\log n)$.
  Indeed,
  \begin{align*}
    \mathbf{a}
    &=
    (y_1,\ldots,y_{n-1},0)'
    -
    \mathbfit{U}_n
    \left[
      (\mathbfit{U}'_n\mathbfit{U}_n+a\mathbfit{I})^{-1}
      \mathbfit{U}_n'
      (y_1,\ldots,y_{n-1},0)'
    \right]\\
  \intertext{and}
    \mathbf{b}
    &=
    (0,\ldots,0,1)'
    -
    \mathbfit{U}_n
    \left[
      (\mathbfit{U}'_n\mathbfit{U}_n+a\mathbfit{I})^{-1}
      \mathbfit{U}_n'
      (0,\ldots,0,1)'
    \right]
  \end{align*}
  can be computed in time $O(n)$,
  and sorting $P$ can be done in time $O(n\log n)$
  (see, e.g., \cite{cormen/etal:2001}, Part~II).

  \ifarXiv
    \section*{Appendix~C: Instruction manual for the R package \texttt{PredictiveRegression}}
  \fi

  The R package \texttt{PredictiveRegression}
  implements the three prediction algorithms
  (IID predictor, Gauss predictor, and MVA predictor)
  described in this paper.
  It is available from CRAN via {\tt http://www.r-project.org/}.

  The program implementing the IID predictor
  is \texttt{iidpred}.
  It is called using the command
  \texttt{iidpred(train,test,epsilons,ridge)}
  with the default values \texttt{c(0.05,0.01)} for \texttt{epsilons}
  and \texttt{0} for \texttt{ridge}.
  The arguments are:
  \begin{description}
  \item[\texttt{train}]
    The training set as a matrix of size $N\times(K+1)$.
    Each row describes an observation.
    Columns $1$ to $K$ are the explanatory variables,
    and column $K+1$ is the response variables.
  \item[\texttt{test}]
    The test set as a matrix of size $N_2\times K$.
    Each row corresponds to an observation
    (but without the response variable).
    Columns $1$ to $K$ are the explanatory variables.
  \item[\texttt{epsilons}]
    A vector of several significance levels.
    Each significance level $\texttt{epsilons}[j]$ is a number between $0$ and $1$.
    The default value is $(5\%,1\%)'$.
  \item[\texttt{ridge}]
    The ridge coefficient, a nonnegative number.
    The default value is $0$;
    setting it to a small positive constant might lead to more stable results.
  \end{description}
  The output is a list of three elements:
  \begin{description}
  \item[\texttt{output[[1]]}]
    This is the matrix of lower bounds of prediction intervals.
    Its size is $N_2\times N_{\epsilon}$,
    where $N_2$ is the number of test observations
    and $N_{\epsilon}$ is the number of significance levels.
    The element $\texttt{output[[1]]}[i,j]$ of $\texttt{output[[1]]}$
    is the lower bound $a$ of the prediction interval $[a,b]$
    for the $i$th test observation
    and for the $j$th significance level $\texttt{epsilons}[j]$
    in the vector \texttt{epsilons}.
  \item[\texttt{output[[2]]}]
    The matrix of upper bounds $b$,
    with the same structure as \texttt{output[[1]]}.
    Typically $a = \texttt{output[[1]]}[i,j]$
    and $b = \texttt{output[[2]]}[i,j]$ are real numbers
    such that $a\le b$.
    Exceptions: $a$ is allowed to be $-\infty$
    and $b$ is allowed to be $\infty$;
    the only case where $a > b$
    is $a=\infty$ and $b=-\infty$
    (the empty prediction $[a,b]$).
  \item[\texttt{output[[3]]}]
    The termination code, which is one of:
    \begin{description}
    \item[0] normal termination;
    \item[1] illegal parameters
      (the training and test sets have different numbers of explanatory variables);
    \item[2] too few observations for all significance levels.
    \end{description}
  \end{description}
  This program implements the algorithm described on pp.~30--34
  of \cite{vovk/etal:2005book}
  and, in more detail, in the previous section.
  A simple example of its use is
  \begin{verbatim}
    train <- matrix(c(0,10,20,30, 1.01,10.99,21.01,30.99),
       nrow=4, ncol=2);
    test <- matrix(c(5,15,25), nrow=3, ncol=1);
    output <- iidpred(train,test,c(0.05,0.2),0.01);
    print(output[[1]]);
    print(output[[2]]);
  \end{verbatim}

  The programs implementing the MVA predictor and the Gauss predictor
  are \texttt{mvapred} and \texttt{gausspred}, respectively.
  They are called in the same way as \texttt{iidpred},
  with the only exception
  that \texttt{gausspred} does not require the parameter \texttt{ridge}.
  The termination code for \texttt{mvapred} and \texttt{gausspred}
  (returned in \texttt{output[[3]]}) is:
  \begin{description}
  \item[0] normal termination;
  \item[1] illegal parameters
    (the training and test sets have different numbers of explanatory variables);
  \item[2] too few observations.
  \end{description}
  An example of use of \texttt{mvapred} is
  \begin{verbatim}
    train <- matrix(c(0,10,20,30, 1.01,10.99,21.01,30.99),
       nrow=4,ncol=2);
    test <- matrix(c(5,15,25), nrow=3, ncol=1);
    output <- mvapred(train,test,c(0.05,0.2),0.01);
    print(output[[1]]);
    print(output[[2]]);
  \end{verbatim}
  and an example of use of \texttt{gausspred} is
  \begin{verbatim}
    train <- matrix(c(1,2,3,4, 2.01,2.99,4.01,4.99),
       nrow=4, ncol=2);
    test <- matrix(c(0,10,20), nrow=3, ncol=1);
    output <- gausspred(train,test,c(0.05,0.2));
    print(output[[1]]);
    print(output[[2]]);
  \end{verbatim}

  \ifarXiv
    \section*{Appendix~D: Goodness of definitions of nonconformity measures}
  \fi

  In this appendix we will check in detail that all our nonconformity measures
  $A_n$ depend on the first $n-1$ observations only through the value of $S_{n-1}$
  (briefly this was explained in the main part of the paper).
  For the IID model and the IID--Gauss model
  this is obvious.

  In the case of the Gauss linear model,
  to check that the expression on the right-hand side of (\ref{eq:nonconformityGLM})
  depends on the first $n-1$ observations only through the value of $S_{n-1}$,
  we can first check this for $\mathbfit{Z}'_{n-1}\mathbfit{Z}_{n-1}$,
  then for $\hat{\boldsymbol{\gamma}}_{n-1}$,
  and finally for $\hat{y}_n$ and $\hat{\sigma}^2_{n-1}$
  (the definition of the latter has to be expanded).

  In the case of the MVA model,
  we are required to check that the right-hand side of (\ref{eq:nonconformityMVA})
  depends on the first $n-1$ observations only through the value of $S_{n-1}$.
  According to (\ref{eq:residuals2}),
  $\mathbf{e}=\mathbf{y}-\mathbfit{U}\mathbf{c}$,
  where $\mathbf{c}$ is a known vector.
  First notice that $e_n$ depends only on the $n$th observation
  and that $\overline{e}_{n-1}$ depends only on the average of $y_1,\ldots,y_{n-1}$
  and the average of $\mathbf{x}_1,\ldots,\mathbf{x}_{n-1}$.
  Therefore, both $e_n$ and $\overline{e}_{n-1}$
  (equivalently, $\sum_{i=1}^{n-1}e_i$)
  depend on the first $n-1$ observations only through the value of $S_{n-1}$.
  It remains to show that $\sum_{i=1}^{n-1}e_i^2$
  depends on the first $n-1$ observations only through the value of $S_{n-1}$.
  This immediately follows from the definition of $S_{n-1}$.

\end{document}